\catcode`!=11 

\newcount\aux%
\newcount\auxa%
\newcount\auxb%
\newcount\m%
\newcount\n%
\newcount\x%
\newcount\y%
\newcount\xl%
\newcount\yl%
\newcount\d%
\newcount\dnm%
\newcount\xa%
\newcount\xb%
\newcount\xmed%
\newcount\xc%
\newcount\xd%
\newcount\ya%
\newcount\yb%
\newcount\ymed%
\newcount\yc%
\newcount\yd
\newcount\expansao%
\newcount\tipografo
\newcount\distanciaobjmor
\newcount\tipoarco
\newif\ifpara%
\newbox\caixa%
\newbox\caixaaux%
\newif\ifnvazia%
\newif\ifvazia%
\newif\ifcompara%
\newif\ifdiferentes%
\newcount\xaux%
\newcount\yaux%
\newcount\guardaauxa%
\newcount\alt%
\newcount\larg%
\newcount\prof%
\newcount\auxqx
\newcount\auxqy
\newif\ifajusta%
\newif\ifajustadist
\def\objPartida{}%
\def\objChegada{}%
\def\objNulo{}%


\def\!vazia{:}

\def\!pilhanvazia#1{\let\arg=#1%
\if:\arg\ \nvaziafalse\vaziatrue \else \nvaziatrue\vaziafalse\fi}

\def\!coloca#1#2{\edef\pilha{#1.#2}}

\def\!guarda(#1)(#2,#3)(#4,#5,#6){\def\id{#1}%
\xaux=#2%
\yaux=#3%
\alt=#4%
\larg=#5%
\prof=#6%
}

\def\!topaux#1.#2:{\!guarda#1}
\def\!topo#1{\expandafter\!topaux#1}

\def\!popaux#1.#2:{\def\pilha{#2:}}
\def\!retira#1{\expandafter\!popaux#1}

\def\!comparaaux#1#2{\let\argA=#1\let\argB=#2%
\ifx\argA\argB\comparatrue\diferentesfalse\else\comparafalse\diferentestrue\fi}

\def\!compara#1#2{\!comparaaux{#1}{#2}}

\def\!absoluto#1#2{\n=#1%
  \ifnum \n > 0
    #2=\n
  \else
    \multiply \n by -1
    #2=\n
  \fi}






\def\!ajusta#1#2#3#4#5#6{\aux=#5%
  \let\auxobj=#6%
  \ifcase \tipografo    
    \ifnum\number\aux=10 
      \ajustadisttrue 
    \else
      \ajustadistfalse  
    \fi
  \else  
   \ajustadistfalse
  \fi
  \ifajustadist
   \let\pilhaaux=\pilha%
   \loop%
     \!topo{\pilha}%
     \!retira{\pilha}%
     \!compara{\id}{\auxobj}%
     \ifcompara\nvaziafalse \else\!pilhanvazia\pilha \fi%
     \ifnvazia%
   \repeat%
   \let\pilha=\pilhaaux%
   \ifvazia%
    \ifdiferentes%
     \larg=1310720
     \prof=655360%
     \alt=655360%
    \fi%
   \fi%
   \divide\larg by 131072
   \divide\prof by 65536
   \divide\alt by 65536
   \ifnum\number\y=\number\yl
    \advance\larg by 3
    \ifnum\number\larg>\aux
     #5=\larg
    \fi
   \else
    \ifnum\number\x=\number\xl
     \ifnum\number\yl>\number\y
      \ifnum\number\alt>\aux
       #5=\alt
      \fi
     \else
      \advance\prof by 5
      \ifnum\number\prof>\aux
       #5=\prof
      \fi
     \fi
    \else
     \auxqx=\x
     \advance\auxqx by -\xl
     \!absoluto{\auxqx}{\auxqx}%
     \auxqy=\y
     \advance\auxqy by -\yl
     \!absoluto{\auxqy}{\auxqy}%
     \ifnum\auxqx>\auxqy
      \ifnum\larg<10
       \larg=10
      \fi
      \advance\larg by 3
      #5=\larg
     \else
      \ifnum\yl>\y
       \ifnum\larg<10
        \larg=10
       \fi
      \advance\alt by 6
       #5=\alt
      \else
      \advance\prof by 11
       #5=\prof
      \fi
     \fi
    \fi
   \fi
\fi} 

\def\!raiz#1#2{\n=#1%
  \m=1%
  \loop
    \aux=\m%
    \advance \aux by 1%
    \multiply \aux by \aux%
    \ifnum \aux < \n%
      \advance \m by 1%
      \paratrue%
    \else\ifnum \aux=\n%
      \advance \m by 1%
      \paratrue%
       \else\parafalse%
       \fi
    \fi
  \ifpara%
  \repeat
#2=\m}

\def\!ucoord#1#2#3#4#5#6#7{\aux=#2%
  \advance \aux by -#1%
  \multiply \aux by #4%
  \divide \aux by #5%
  \ifnum #7 = -1 \multiply \aux by -1 \fi%
  \advance \aux by #3%
#6=\aux}

\def\!quadrado#1#2#3{\aux=#1%
  \advance \aux by -#2%
  \multiply \aux by \aux%
#3=\aux}

\def\!distnomemor#1#2#3#4#5#6{\setbox0=\hbox{#5}%
  \aux=#1
  \advance \aux by -#3
  \ifnum \aux=0
     \aux=\wd0 \divide \aux by 131072
     \advance \aux by 3
     #6=\aux
  \else
     \aux=#2
     \advance \aux by -#4
     \ifnum \aux=0
        \aux=\ht0 \advance \aux by \dp0 \divide \aux by 131072
        \advance \aux by 3
        #6=\aux%
     \else
     #6=3
     \fi
   \fi
}

\def\begindc#1{\!ifnextchar[{\!begindc{#1}}{\!begindc{#1}[30]}}
\def\!begindc#1[#2]{\beginpicture 
  \let\pilha=\!vazia
  \setcoordinatesystem units <1pt,1pt>
  \expansao=#2
  \ifcase #1
    \distanciaobjmor=10
    \tipoarco=0         
    \tipografo=0        
  \or
    \distanciaobjmor=2
    \tipoarco=0         
    \tipografo=1        
  \or
    \distanciaobjmor=1
    \tipoarco=2         
    \tipografo=2        
  \or
    \distanciaobjmor=8
    \tipoarco=0         
    \tipografo=3        
  \or
    \distanciaobjmor=8
    \tipoarco=2         
    \tipografo=4        
  \fi}

\def\enddc{\endpicture}

\def\mor{%
  \!ifnextchar({\!morxy}{\!morObjA}}
\def\!morxy(#1,#2){%
  \!ifnextchar({\!morxyl{#1}{#2}}{\!morObjB{#1}{#2}}}
\def\!morxyl#1#2(#3,#4){%
  \!ifnextchar[{\!mora{#1}{#2}{#3}{#4}}{\!mora{#1}{#2}{#3}{#4}[\number\distanciaobjmor,\number\distanciaobjmor]}}%
\def\!morObjA#1{%
 \let\pilhaaux=\pilha%
 \def\objPartida{#1}%
 \loop%
    \!topo\pilha%
    \!retira\pilha%
    \!compara{\id}{\objPartida}%
    \ifcompara \nvaziafalse \else \!pilhanvazia\pilha \fi%
   \ifnvazia%
 \repeat%
 \ifvazia%
  \ifdiferentes%
   Error: Incorrect label specification%
   \xaux=1%
   \yaux=1%
  \fi%
 \fi%
 \let\pilha=\pilhaaux%
 \!ifnextchar({\!morxyl{\number\xaux}{\number\yaux}}{\!morObjB{\number\xaux}{\number\yaux}}}
\def\!morObjB#1#2#3{%
  \x=#1
  \y=#2
 \def\objChegada{#3}%
 \let\pilhaaux=\pilha%
 \loop
    \!topo\pilha %
    \!retira\pilha%
    \!compara{\id}{\objChegada}%
    \ifcompara \nvaziafalse \else \!pilhanvazia\pilha \fi
   \ifnvazia
 \repeat
 \ifvazia
  \ifdiferentes%
   Error: Incorrect label specification
   \xaux=\x%
   \advance\xaux by \x%
   \yaux=\y%
   \advance\yaux by \y%
  \fi
 \fi
 \let\pilha=\pilhaaux
 \!ifnextchar[{\!mora{\number\x}{\number\y}{\number\xaux}{\number\yaux}}{\!mora{\number\x}{\number\y}{\number\xaux}{\number\yaux}[\number\distanciaobjmor,\number\distanciaobjmor]}}
\def\!mora#1#2#3#4[#5,#6]#7{%
  \!ifnextchar[{\!morb{#1}{#2}{#3}{#4}{#5}{#6}{#7}}{\!morb{#1}{#2}{#3}{#4}{#5}{#6}{#7}[1,\number\tipoarco] }}
\def\!morb#1#2#3#4#5#6#7[#8,#9]{\x=#1%
  \y=#2%
  \xl=#3%
  \yl=#4%
  \multiply \x by \expansao%
  \multiply \y by \expansao%
  \multiply \xl by \expansao%
  \multiply \yl by \expansao%
  \!quadrado{\number\x}{\number\xl}{\auxa}%
  \!quadrado{\number\y}{\number\yl}{\auxb}%
  \d=\auxa%
  \advance \d by \auxb%
  \!raiz{\d}{\d}%
  \auxa=#5
  \!compara{\objNulo}{\objPartida}%
  \ifdiferentes
   \!ajusta{\x}{\xl}{\y}{\yl}{\auxa}{\objPartida}%
   \ajustatrue
   \def\objPartida{}
  \fi
  \guardaauxa=\auxa
  \!ucoord{\number\x}{\number\xl}{\number\x}{\auxa}{\number\d}{\xa}{1}%
  \!ucoord{\number\y}{\number\yl}{\number\y}{\auxa}{\number\d}{\ya}{1}%
  \auxa=\d%
  \auxb=#6
  \!compara{\objNulo}{\objChegada}%
  \ifdiferentes
   \!ajusta{\x}{\xl}{\y}{\yl}{\auxb}{\objChegada}%
   \def\objChegada{}
  \fi
  \advance \auxa by -\auxb%
  \!ucoord{\number\x}{\number\xl}{\number\x}{\number\auxa}{\number\d}{\xb}{1}%
  \!ucoord{\number\y}{\number\yl}{\number\y}{\number\auxa}{\number\d}{\yb}{1}%
  \xmed=\xa%
  \advance \xmed by \xb%
  \divide \xmed by 2
  \ymed=\ya%
  \advance \ymed by \yb%
  \divide \ymed by 2
  \!distnomemor{\number\x}{\number\y}{\number\xl}{\number\yl}{#7}{\dnm}%
  \!ucoord{\number\y}{\number\yl}{\number\xmed}{\number\dnm}{\number\d}{\xc}{-#8}%
  \!ucoord{\number\x}{\number\xl}{\number\ymed}{\number\dnm}{\number\d}{\yc}{#8}%
\ifcase #9  
  \arrow <4pt> [.2,1.1] from {\xa} {\ya} to {\xb} {\yb}
\or  
  \setdashes
  \arrow <4pt> [.2,1.1] from {\xa} {\ya} to {\xb} {\yb}
  \setsolid
\or  
  \setlinear
  \plot {\xa} {\ya}  {\xb} {\yb} /
\or  
  \auxa=\guardaauxa
  \advance \auxa by 3%
 \!ucoord{\number\x}{\number\xl}{\number\x}{\number\auxa}{\number\d}{\xa}{1}%
 \!ucoord{\number\y}{\number\yl}{\number\y}{\number\auxa}{\number\d}{\ya}{1}%
 \!ucoord{\number\y}{\number\yl}{\number\xa}{3}{\number\d}{\xd}{-1}%
 \!ucoord{\number\x}{\number\xl}{\number\ya}{3}{\number\d}{\yd}{1}%
  \arrow <4pt> [.2,1.1] from {\xa} {\ya} to {\xb} {\yb}
  \circulararc -180 degrees from {\xa} {\ya} center at {\xd} {\yd}
\or  
  \auxa=3
 \!ucoord{\number\y}{\number\yl}{\number\xa}{\number\auxa}{\number\d}{\xmed}{-1}%
 \!ucoord{\number\x}{\number\xl}{\number\ya}{\number\auxa}{\number\d}{\ymed}{1}%
 \!ucoord{\number\y}{\number\yl}{\number\xa}{\number\auxa}{\number\d}{\xd}{1}%
 \!ucoord{\number\x}{\number\xl}{\number\ya}{\number\auxa}{\number\d}{\yd}{-1}%
  \arrow <4pt> [.2,1.1] from {\xa} {\ya} to {\xb} {\yb}
  \setlinear
  \plot {\xmed} {\ymed}  {\xd} {\yd} /
\fi
\auxa=\xl
\advance \auxa by -\x%
\ifnum \auxa=0 
  \put {#7} at {\xc} {\yc}
\else
  \auxb=\yl
  \advance \auxb by -\y%
  \ifnum \auxb=0 \put {#7} at {\xc} {\yc}
  \else 
    \ifnum \auxa > 0 
      \ifnum \auxb > 0
        \ifnum #8=1
          \put {#7} [rb] at {\xc} {\yc}
        \else 
          \put {#7} [lt] at {\xc} {\yc}
        \fi
      \else
        \ifnum #8=1
          \put {#7} [lb] at {\xc} {\yc}
        \else 
          \put {#7} [rt] at {\xc} {\yc}
        \fi
      \fi
    \else
      \ifnum \auxb > 0 
        \ifnum #8=1
          \put {#7} [rt] at {\xc} {\yc}
        \else 
          \put {#7} [lb] at {\xc} {\yc}
        \fi
      \else
        \ifnum #8=1
          \put {#7} [lt] at {\xc} {\yc}
        \else 
          \put {#7} [rb] at {\xc} {\yc}
        \fi
      \fi
    \fi
  \fi
\fi
}

\def\modifplot(#1{\!modifqcurve #1}
\def\!modifqcurve(#1,#2){\x=#1%
  \y=#2%
  \multiply \x by \expansao%
  \multiply \y by \expansao%
  \!start (\x,\y)
  \!modifQjoin}
\def\!modifQjoin(#1,#2)(#3,#4){\x=#1%
  \y=#2%
  \xl=#3%
  \yl=#4%
  \multiply \x by \expansao%
  \multiply \y by \expansao%
  \multiply \xl by \expansao%
  \multiply \yl by \expansao%
  \!qjoin (\x,\y) (\xl,\yl)             
  \!ifnextchar){\!fim}{\!modifQjoin}}
\def\!fim){\ignorespaces}

\def\setaxy(#1{\!pontosxy #1}
\def\!pontosxy(#1,#2){%
  \!maispontosxy}
\def\!maispontosxy(#1,#2)(#3,#4){%
  \!ifnextchar){\!fimxy#3,#4}{\!maispontosxy}}
\def\!fimxy#1,#2){\x=#1%
  \y=#2
  \multiply \x by \expansao
  \multiply \y by \expansao
  \xl=\x%
  \yl=\y%
  \aux=1%
  \multiply \aux by \auxa%
  \advance\xl by \aux%
  \aux=1%
  \multiply \aux by \auxb%
  \advance\yl by \aux%
  \arrow <4pt> [.2,1.1] from {\x} {\y} to {\xl} {\yl}}

\def\cmor#1 #2(#3,#4)#5{%
  \!ifnextchar[{\!cmora{#1}{#2}{#3}{#4}{#5}}{\!cmora{#1}{#2}{#3}{#4}{#5}[0] }}
\def\!cmora#1#2#3#4#5[#6]{%
  \ifcase #2
      \auxa=0
      \auxb=1
    \or
      \auxa=0
      \auxb=-1
    \or
      \auxa=1
      \auxb=0
    \or
      \auxa=-1
      \auxb=0
    \fi  
  \ifcase #6  
    \modifplot#1
    \setaxy#1
  \or  
    \setdashes
    \modifplot#1
    \setaxy#1
    \setsolid
  \or  
    \modifplot#1
  \fi  
  \x=#3%
  \y=#4%
  \multiply \x by \expansao%
  \multiply \y by \expansao%
  \put {#5} at {\x} {\y}}

\def\obj(#1,#2){%
  \!ifnextchar[{\!obja{#1}{#2}}{\!obja{#1}{#2}[Nulo]}}
\def\!obja#1#2[#3]#4{%
  \!ifnextchar[{\!objb{#1}{#2}{#3}{#4}}{\!objb{#1}{#2}{#3}{#4}[1]}}
\def\!objb#1#2#3#4[#5]{%
  \x=#1%
  \y=#2%
  \def\!pinta{\normalsize$\bullet$}
  \def\!nulo{Nulo}%
  \def\!arg{#3}%
  \!compara{\!arg}{\!nulo}%
  \ifcompara\def\!arg{#4}\fi%
  \multiply \x by \expansao%
  \multiply \y by \expansao%
  \setbox\caixa=\hbox{#4}%
  \!coloca{(\!arg)(#1,#2)(\number\ht\caixa,\number\wd\caixa,\number\dp\caixa)}{\pilha}%
  \auxa=\wd\caixa \divide \auxa by 131072 
  \advance \auxa by 5
  \auxb=\ht\caixa
  \advance \auxb by \number\dp\caixa
  \divide \auxb by 131072 
  \advance \auxb by 5
  \ifcase \tipografo    
    \put{#4} at {\x} {\y}
  \or                   
    \ifcase #5 
      \put{#4} at {\x} {\y}
    \or        
      \put{\!pinta} at {\x} {\y}
      \advance \y by \number\auxb  
      \put{#4} at {\x} {\y}
    \or        
      \put{\!pinta} at {\x} {\y}
      \advance \auxa by -2  
      \advance \auxb by -2  
      \advance \x by \number\auxa  
      \advance \y by \number\auxb  
      \put{#4} at {\x} {\y}   
    \or        
      \put{\!pinta} at {\x} {\y}
      \advance \x by \number\auxa  
      \put{#4} at {\x} {\y}   
    \or        
      \put{\!pinta} at {\x} {\y}
      \advance \auxa by -2  
      \advance \auxb by -2  
      \advance \x by \number\auxa  
      \advance \y by -\number\auxb  
      \put{#4} at {\x} {\y}   
    \or        
      \put{\!pinta} at {\x} {\y}
      \advance \y by -\number\auxb  
      \put{#4} at {\x} {\y}   
    \or        
      \put{\!pinta} at {\x} {\y}
      \advance \auxa by -2  
      \advance \auxb by -2  
      \advance \x by -\number\auxa  
      \advance \y by -\number\auxb  
      \put{#4} at {\x} {\y}   
    \or        
      \put{\!pinta} at {\x} {\y}
      \advance \x by -\number\auxa  
      \put{#4} at {\x} {\y}   
    \or        
      \put{\!pinta} at {\x} {\y}
      \advance \auxa by -2  
      \advance \auxb by -2  
      \advance \x by -\number\auxa  
      \advance \y by \number\auxb  
      \put{#4} at {\x} {\y}   
    \fi
  \or                   
    \ifcase #5 
      \put{#4} at {\x} {\y}
    \or        
      \put{\!pinta} at {\x} {\y}
      \advance \y by \number\auxb  
      \put{#4} at {\x} {\y}
    \or        
      \put{\!pinta} at {\x} {\y}
      \advance \auxa by -2  
      \advance \auxb by -2  
      \advance \x by \number\auxa  
      \advance \y by \number\auxb  
      \put{#4} at {\x} {\y}   
    \or        
      \put{\!pinta} at {\x} {\y}
      \advance \x by \number\auxa  
      \put{#4} at {\x} {\y}   
    \or        
      \put{\!pinta} at {\x} {\y}
      \advance \auxa by -2  
      \advance \auxb by -2
      \advance \x by \number\auxa  
      \advance \y by -\number\auxb 
      \put{#4} at {\x} {\y}   
    \or        
      \put{\!pinta} at {\x} {\y}
      \advance \y by -\number\auxb 
      \put{#4} at {\x} {\y}   
    \or        
      \put{\!pinta} at {\x} {\y}
      \advance \auxa by -2  
      \advance \auxb by -2
      \advance \x by -\number\auxa 
      \advance \y by -\number\auxb 
      \put{#4} at {\x} {\y}   
    \or        
      \put{\!pinta} at {\x} {\y}
      \advance \x by -\number\auxa 
      \put{#4} at {\x} {\y}   
    \or        
      \put{\!pinta} at {\x} {\y}
      \advance \auxa by -2  
      \advance \auxb by -2
      \advance \x by -\number\auxa 
      \advance \y by \number\auxb  
      \put{#4} at {\x} {\y}   
    \fi
   \else 
     \ifnum\auxa<\auxb 
       \aux=\auxb
     \else
       \aux=\auxa
     \fi
     \ifdim\wd\caixa<1em
       \dimen99 = 1em
       \aux=\dimen99 \divide \aux by 131072 
       \advance \aux by 5
     \fi
     \advance\aux by -2 
     \multiply\aux by 2 %
     \ifnum\aux<30
       \put{\circle{\aux}} [Bl] at {\x} {\y}
     \else
       \multiply\auxa by 2
       \multiply\auxb by 2
       \put{\oval(\auxa,\auxb)} [Bl] at {\x} {\y}
     \fi
     \put{#4} at {\x} {\y}
   \fi   
}

\catcode`!=12 
 
\tolerance=10000
\hsize=15truecm \hoffset 0.4truecm           
\vsize=21truecm \voffset 1.5truecm


\def\Abstract#1\par{\centerline{\vbox{\hsize 11truecm
    \noindent {\bf Abstract.} #1\hfil}} \vskip 0.75cm}

\newcount\numsezione
\def\Sezione#1.{\advance\numsezione by1
    \vskip1cm\goodbreak\noindent
    {\bf \the\numsezione\ - #1. \hfill}
    \vskip 0.25cm \nobreak}

\def\References{\vskip1cm\goodbreak\noindent
    {\bf References \hfill}
    \vskip 0.25cm \nobreak}

\long\def\proclaim#1. #2\par{\goodbreak\noindent
    {\bf #1. }
    {\sl#2}
    \par\ifdim \lastskip <\medskipamount \removelastskip \penalty55 \bigskip \fi}

\def\bull{{\vrule height.9ex width.8ex depth-.1ex}}
\def\s{\scriptscriptstyle}

\def\b#1{{\bf #1}}
\def\lra{\longrightarrow}

\font\tenmsbm=msbm10
\font\sevenmsbm=msbm7
\font\fivemsbm=msbm5
\newfam\amsBfam
\textfont\amsBfam=\tenmsbm
\scriptfont\amsBfam=\sevenmsbm
\scriptscriptfont\amsBfam=\fivemsbm
\def\bbb{\fam\the\amsBfam\tenmsbm}

\font\tenmsam=msam10
\font\sevenmsam=msam7
\font\fivemsam=msam5
\newfam\amsAfam
\textfont\amsAfam=\tenmsam
\scriptfont\amsAfam=\sevenmsam
\scriptscriptfont\amsAfam=\fivemsam
\def\aaa{\fam\the\amsAfam\tenmsam}

\mathchardef\gul"3\the\amsAfam52
\mathchardef\lug"3\the\amsAfam51
\mathchardef\R"5\the\amsBfam52

\def\R {\bbb{R}}
\def\N {\bbb{N}}

\def\zd {\delta}
\def\zD {\Delta}
\def\zl {\lambda}
\def\zq {\chi}

\def\rd {\rm d}
\def\rD {\rm D}

		\catcode`\"=12
	\font\kropa=lcircle10 scaled 1700
	\def\ybl{\setbox0=\hbox{\kropa \char"70} \kern1.5pt \raise.35pt \box0}
		\catcode`\"=\active

	\def\tdot{\ybl}

\input pictex

\font\bfsymbtext=cmmib10 
\font\bfsymbscript=cmmib10 at 7pt
\font\bfsymbscriptscript=cmmib10 at 5 pt
\def\boldsymbol#1{{\mathchoice%
	{\hbox{$\displaystyle\textfont1=\bfsymbtext\scriptfont1=%
		\bfsymbscript\scriptscriptfont1=\bfsymbscriptscript #1$}}%
	{\hbox{$\textstyle\textfont1=\bfsymbtext\scriptfont1=%
		\bfsymbscript\scriptscriptfont1=\bfsymbscriptscript #1$}}%
	{\hbox{$\scriptstyle\scriptfont1=\bfsymbscript%
		\scriptscriptfont1=\bfsymbscriptscript #1$}}%
	{\hbox{$\scriptscriptstyle%
	\scriptscriptfont1=\bfsymbscriptscript {#1}$}}}}

\def\bu{\boldsymbol{\upsilon}}


\centerline {\bf{Polarizations and differential calculus in affine spaces.}}
\vskip 2\baselineskip 
\centerline {MARGHERITA BARILE$^1$}
\centerline{FIORELLA BARONE$^2$}
\centerline{WLODZIMIERZ M. TULCZYJEW$^3$\footnote{(*)}{Supported by PRIN SINTESI.}}
\centerline{\hphantom{X}}
\centerline {$^1$Dipartimento di Matematica, Via Orabona, 4}
\centerline {Universit\`a di Bari, 70125 Bari, Italy}
\centerline {barile@dm.uniba.it}
\centerline {Tel. 0039 0805442711 -- Fax 0039 0805963612} 
\centerline{\hphantom{X}}
\centerline {$^2$Dipartimento di Matematica, Via Orabona, 4}
\centerline {Universit\`a di Bari, 70125 Bari, Italy}
\centerline {barone@dm.uniba.it} 
\centerline {Tel. 0039 0805442711 -- Fax 0039 0805963612}
\centerline{\hphantom{X}}
\centerline {$^3$associated with} 
\centerline {Istituto Nazionale di Fisica Nucleare, Sezione di Napoli, Italy}
\centerline {Universit\`a di Camerino, Italy}
\centerline{Valle San Benedetto, 2}
\centerline{62030 Monte Cavallo (MC), Italy}
\centerline {tulczy@libero.it} 
\centerline {Tel. 0039 0737519748 -- Fax 0039 0737402529}

\vskip 3\baselineskip
\Abstract 
Within the framework of mappings between affine spaces, the notion of $n$-th polarization of a function will lead to an intrinsic  characterization of polynomial functions.
We prove that  the characteristic features of derivations, such as linearity, iterability, Leibniz and chain rules, are shared -- at the finite level -- by the polarization operators. 
We give these results by means of explicit general formulae, which are valid at any order $n$, and are based on combinatorial identities.
The infinitesimal limits of the $n$-th polarizations of a function will yield its $n$-th derivatives (without resorting to the usual recursive definition), and the above mentioned properties will be recovered directly in the limit.  
Polynomial functions will allow us to produce a coordinate free version of Taylor's formula.

\vskip 2\baselineskip

\vskip 1\baselineskip
\noindent A.M.S. CLASSIFICATION: 26C05, 05A10, 26B12, 47H60.

\vskip 3\baselineskip
\goodbreak

\Sezione Introduction.

As it is well known, the notion of increment (or finite difference) of a  function $f$ is based on the affine structure of the space where $f$ is defined.
Using the affine structure, we generalize the idea (Sec.2) and directly define the $n$-th increment of $f$, starting at a point $q$, in $n$ arbitrary directions. 
This is the value of the $n$-th {\it polarization} of $f$ at $q$ in those directions.

The resulting polarization operators will satisfy (Sec.3) all the main characteristic properties of multidirectional derivations, i.e., linearity, iterability, Leibniz and chain rules.

A more refined study  will then show how affine spaces are the natural framework for  both  {\it homogeneous functions} and  {\it homogeneous polynomials} --  the latter being, essentially,  polarizations evaluated  along diagonals, i.e.,  increments in only one direction. 

For homogeneous functions an analogue of the Euler theorem is proved (Sec.4).

For homogeneous polynomials a characterization theorem is proved (Sec.5).

Sums of homogeneous polynomials will finally define {\it polynomial functions} on an affine space.

A natural application of the above algebraic study is the intrinsic differential calculus on affine spaces.
The limit of the $n$-th polarization will be (Sec.6) the $n$-th multidirectional derivative, defined directly without iterations. 
We retrieve the usual properties of the $n$-th derivatives as the infinitesimal limits of the corresponding properties of  the $n$-th polarizations (Sec.7 -- 10). 

The role of the polynomial functions will be played in the final section (Sec.11) where differentiability is treated. 
Polynomials will be the tool for an intrinsic Taylor formula. 

This paper is the starting point for an algebraic approach to jet spaces (both Weil and Ehresmann type) whose differential structures will be modelled on affine spaces.

\vskip 2\baselineskip
\centerline {\bf{A. The calculus of polarizations in affine spaces.}}

\Sezione Polarization of functions on affine spaces.

Let $Q$ be an affine space modelled on a vector space $V$.  For each $n \in \N$ we consider the vector space $B^n(Q)$ of real functions on the product $Q\times V^n$.  The space $A(Q) = B^0(Q)$ is a commutative, associative algebra. It is the algebra of all functions on $Q$.

Given $F_1 \in B^{n_1}(Q)$ and $F_2 \in B^{n_2}(Q)$ we define the product $F_1F_2 \in B^{n_1+n_2}(Q)$ by
		$$F_1F_2(q; v_1,\dots, v_{n_1}, v_{n_1+1}, \dots, v_{n_1+n_2}) = F_1(q; v_1,\dots, v_{n_1}) F_2(q; v_{n_1+1}, \dots,
v_{n_1+n_2}) \eqno(2)$$
This definition extends the product in $B^0(Q)$.  

The direct sum
		$$B(Q) = \oplus_{n=0}^{\infty} B^n(Q) \eqno(3)$$
is an associative algebra.  It is not commutative.

A function $F\in B^n(Q)$ is said to be {\it symmetric} if for every permutation $\sigma\in S_n$ the equality
		$$F(q; v_{\sigma(1)},\dots, v_{\sigma(n)}) = F(q; v_1,\dots, v_n) \eqno(4)$$
holds for each $(q; v_1,\dots, v_n)\in Q\times V^n$.

A function $F\in B^n(Q)$ is said to be {\it multilinear} at $q\in Q$ if for each $1 \le k \le n$ the mapping
		$$\eqalign{
F(q; v_1,\dots,v_{k-1}, \,\cdot\, , v_{k+1}, \ldots, v_n) &: V \to  R \cr
                        &: v \; \mapsto \; F(q; v_1,\ldots,v_{k-1}, \,v\, , v_{k+1}, \ldots v_n)
} \eqno(5)$$
	is linear.

We will use finite index sets $I \subset \N \setminus \{\emptyset\}$ of positive integers.  

The symbol $|I|$ denotes the cardinality of the set $I$.  
The set $\{1,\ldots,n\}$ will be denoted by N. 
 For an index set $I = \{i_1,i_2,\ldots,i_m\}$  we will use the symbol ${\b v}^I$ to denote the sequence $(v_{i_1},v_{i_2},\dots,v_{i_m})$ with $i_1 < i_2< \dots <i_m$.

\vskip 1\baselineskip 
\proclaim Definition 1.												
The $0$-{\it th polarization} of a function $f \in A(Q)$ is the function itself.  For $n > 0$, the $n$-{\it
th polarization} of a function $f \in A(Q)$ is the function
$$\delta^n f \colon Q \times V^n \lra \R \eqno(6)$$
defined by
		$$\delta^n f(q;v_1,v_2,\ldots,v_n) = (-1)^n \sum_{\scriptstyle I \subset N} (-1)^{|I|} f\left(q + 
{\sum_{i \in I}}v_i \right). \eqno(7)$$
By a convenient convention, the term in the above sum corresponding to the empty set $I$ is set equal to $f(q)$.
\hfill$\bull$		

Each index set $I \subset \{1,2,\ldots,n\}$ is uniquely represented by the sequence $(i_1,i_2,\ldots,i_m)$ of its elements arranged in the order of increasing value.  The following alternative definition is obtained by using this representation.  The $n$-{\it th polarization} of a function $f \in A(Q)$ is defined by the alternative formula
		$$\delta^n f(q; v_1, \dots, v_n) = (-1)^n \left[ f(q) + \sum_{m=1}^n (-1)^m\; \sum_{1\le i_1<\dots<i_m\le n} f(q + v_{i_1} +\dots + v_{i_m}) \right]. \eqno(8)$$

The polarization is the sum of
		$$\sum_{m=0}^n {n \choose m} = 2^n \eqno(9)$$
terms.  Each term is the value of $f$ at a point obtained by adding to $q$ the sum of $m$ vectors in the set $\{v_1, v_2,\dots, v_n\}$.  This value appears in the sum with the sign factor $(-1)^m$.  There is an overall sign factor $(-1)^n$.

This description is illustrated in the case $n=3$ by the diagram
\vskip15mm

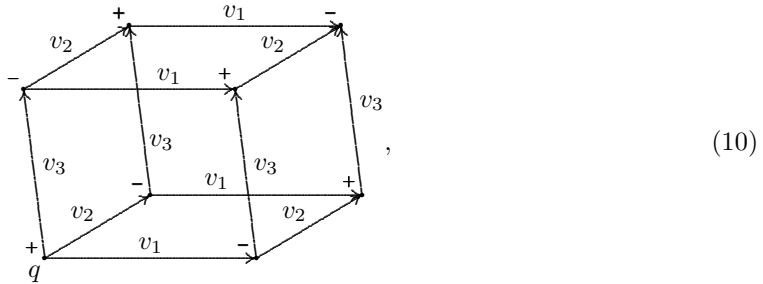
\def\pxb#1{\setbox0\hbox{${#1}$} \copy0 \kern-\wd0 \kern.2pt \copy0 \kern-\wd0 \kern-.2pt \raise.2pt \copy0 \kern-\wd0
\kern.2pt\raise.2pt\box0}

\def\plus{\pxb{\s +}}
\def\minus{\pxb{\s -}}

		$$\vcenter{
	\begindc{0}[1]
	\obj(12,10)[0000]{$q$}[6]
	\obj(60,60)[v]{$v_3$}[0]
	\obj(16,16)[000]{$\tdot$}[6]
	\obj(11,20)[a]{$\plus$}
	\obj(8,80)[001]{$\tdot$}
	\obj(4,84)[b]{$\minus$}
	\obj(96,16)[100]{$\tdot$}
	\obj(91,21)[c]{$\minus$}
	\obj(56,40)[010]{$\tdot$}
	\obj(51,44)[d]{$\minus$}
	\obj(48,104)[011]{$\tdot$}
	\obj(44,109)[e]{$\plus$}
	\obj(88,80)[101]{$\tdot$}
	\obj(84,86)[f]{$\plus$}
	\obj(136,40)[110]{$\tdot$}
	\obj(131,45)[g]{$\plus$}
	\obj(128,104)[111]{$\tdot$}
	\obj(124,109)[h]{$\minus$}
	\mor{000}{100}[0,0]{$v_1$}[2,0]
	\mor{000}{010}[0,0]{$v_2$}[1,0]
	\mor{000}{001}[0,0]{$v_3$}[2,0]
	\mor{010}{110}[0,0]{$\kern-30pt v_1$}[2,0]
	\mor{100}{110}[0,0]{$v_2$}[1,0]
	\mor{100}{101}[0,0]{$v_3$}[2,0]
	\mor{010}{011}[0,0]{}[2,0]
	\mor{110}{111}[0,0]{$v_3$}[2,0]
	\mor{001}{011}[0,0]{$v_2$}[1,0]
	\mor{001}{101}[0,0]{$\kern30pt v_1$}[2,0]
	\mor{101}{111}[0,0]{$v_2$}[1,0]
	\mor{011}{111}[0,0]{$v_1$}[2,0]
	\enddc},
																										\eqno(10)$$
\vskip15mm
\noindent
where we have labeled the points with the minus sign or the plus sign in correspondence with the factor $(-1)^m$.  The polarization corresponding to the diagram is the expression
   $$\eqalign{
	f(q + v_1 + v_2 + v_3) - f(q + v_1 + v_2) - &f(q + v_1 + v_3) - f(q + v_1 + v_2) \cr
					                                      &+ f(q + v_1) + f(q + v_2) + f(q + v_3) - f(q).	
	}\eqno(11)$$

The polarization of a function is defined as a combination of incremented values of the function.  A direct computation yields the formula
  $$\eqalign{
	f(q + v_1 + v_2 + \dots + v_n) &= \sum_{I \subset N} \;\delta^{|I|} f(q;{\b v}^I) \cr
	&= f(q) + \sum_{m=1}^n \;\; \sum_{1\le i_1<\dots<i_m\le n}\; \delta^m f(q;v_{i_1},\dots,v_{i_m})\,,
	}\eqno(12)$$
which inverts the relation between the polarizations of a function and the incremented values of the function by expressing an incremented value of the function as a combination of polarizations.

It follows from the definition that the polarization of a function is symmetric in its vector arguments.

We introduce the $n$-{\it th unidirectional polarization} of a function defined by
	$$\zD^n f \colon Q \times V \rightarrow {\R} \colon (q;v) \mapsto \zd^n f(q;v, v, \dots, v).
	\eqno(13)$$
It is essentially the restriction of the $n$-th polarization to the product of the space $Q$ with the diagonal of $V^n$.

The relation
	$$\zD^n f(q;v) = (-1)^n \sum_{m=0}^n (-1)^m {n\choose m} f(q+mv)
	\eqno(14)$$
holds.

The $n$-th polarization of a function can be considered the result of the application of the {\it n-th polarization operator}
	$$\zd^n : A(Q) \to B^n(Q). \eqno(15)$$
This linear operator has obvious extensions to operators
        $$\zd^n : B^m(Q) \to B^{n+m}(Q) \eqno(16)$$
defined by

	$$\eqalign{
	\zd^n F(q; v_1, \dots,& v_n , w_1, \dots, w_m)\cr 
	& = (-1)^n\sum_{k=0}^n (-1)^k
       \sum_{\scriptstyle I \subset N \atop{\scriptstyle |I| = k}} 
	F\left(q + {\s{\sum_{i \in I}}}v_i; w_1, \dots, w_m\right) \cr
	&= (-1)^n  \sum_{k=0}^n (-1)^k\; \sum_{1\le i_1<\dots<i_k\le n} 
	F(q + v_{i_1} +\dots + v_{i_k} ; w_1, \dots, w_m).
	}\eqno(17)$$
	
In this definition the arguments $(w_1, \dots, w_m)$ of the function $F$ do not participate in the formation of the polarizations.  It follows that the properties of polarizations are shared by the extended operators.  

The extended definition makes repeated applications of polarization operators possible.

\vskip 1\baselineskip 
\proclaim Proposition 1.													
The relation
	$$\zd^{n_1}\,\zd^{n_2} = \zd^{n_1+n_2} \eqno(18)$$
holds for all integers $n_1$ and $n_2$.

\noindent{\it Proof. }
	$$\eqalign{
	\zd^{n_1}\zd^{n_2} & f(q; v^1_1, \dots , v^1_{n_1}, v^2_1, \dots , v^2_{n_2}) \cr
	& = (-1)^{n_1} \sum_{m_1=0}^{n_1} (-1)^{m_1} \sum_{1\le i^1_1<\dots<i^1_{m_1}\le n_1}                        \zd^{n_2}f(q + v^1_{i_1} +\dots + v^1_{i_{m_1}}; v^2_1, \dots , v^2_{n_2}) \cr
 	& = (-1)^{n_1}  \sum_{m_1=0}^{n_1} (-1)^{m_1} \sum_{1\le  i^1_1 < \dots < i^1_{m_1} \le n_1}\cr
	&\phantom{=} \left[  (-1)^{n_2}  \sum_{m_2 = 0}^{n_2} (-1)^{m_2} 
	\sum_{1 \le i^2_1< \dots< i^2_{m_2} \le n_2}
	f(q + v^1_{i_1} +\dots + v^1_{i_{m_1}} + v^2_{i_1} +\dots + v^2_{i_{m_2}})\right] \cr
	& = (-1)^{n_1+ n_2} \sum_{m_1 = 0}^{n_1} \sum_{m_2 = 0}^{n_2}(-1)^{m_1+m_2} \cr
	&\phantom{=}\sum_{1\le i^1_1<\dots<i^1_{m_1}\le n_1}\;\sum_{1\le i^2_1<\dots<i^2_{m_2}\le n_2}
	f(q + v^1_{i_1} +\dots + v^1_{i_{m_1}} + v^2_{i_1} +\dots + v^2_{i_{m_2}}).
	}\eqno(19)$$

On the other hand we have
	$$\eqalign{
	\zd^{n_1+n_2}&f(q; v^1_1, \dots , v^1_{n_1}, v^2_1, \dots , v^2_{n_2})\cr
	& = (-1)^{n_1+n_2} \sum_{m=0}^{n_1+n_2} (-1)^m
	   \sum_{1\le i_1<\dots<i_m\le n_1+n_2}f(q; v_{ i_1}+ \dots + v_{ i_m}).
	}\eqno(20)$$
	
Each term in the sum (19) appears in the combined sum (20) exactly once and with the correct coefficient
$+1$ or $-1$.  We have thus established the equality
	$$\zd^{n_1}\,\zd^{n_2}f(q; v_1, \dots , v_{n_1+n_2}) = \zd^{n_1+n_2} f(q; v_1, \dots , v_{n_1+n_2})
	\eqno(21)$$
for a function $f \in A(Q)$.  This equality extends easily to the equality
	$$\zd^{n_1}\,\zd^{n_2}F(q; v_1, \dots , v_{n_1+n_2}, w_1, \dots , w_m) = 
	\zd^{n_1+n_2} F(q; v_1, \dots , v_{n_1+n_2}, w_1, \dots , w_m) \eqno(22)$$
for a function $F \in B^m(Q)$.
\hfill$\bull$

The relation
	$$\zD^{n_1}\,\zD^{n_2} = \zD^{n_1+n_2} \eqno(23)$$
is an immediate consequence of the above proposition.

The diagram
\vskip15mm

		$$\vcenter{
	\begindc{0}[1]
	\obj(12,10)[0000]{$q$}[6]
	\obj(60,60)[v]{$v_1^1$}[0]
	\obj(16,16)[000]{$\tdot$}[6]
	\obj(11,20)[a]{$\plus$}
	\obj(8,80)[001]{$\tdot$}
	\obj(4,84)[b]{$\minus$}
	\obj(96,16)[100]{$\tdot$}
	\obj(91,21)[c]{$\minus$}
	\obj(56,40)[010]{$\tdot$}
	\obj(51,44)[d]{$\minus$}
	\obj(48,104)[011]{$\tdot$}
	\obj(44,109)[e]{$\plus$}
	\obj(88,80)[101]{$\tdot$}
	\obj(84,86)[f]{$\plus$}
	\obj(136,40)[110]{$\tdot$}
	\obj(131,45)[g]{$\plus$}
	\obj(128,104)[111]{$\tdot$}
	\obj(124,109)[h]{$\minus$}
	\mor{000}{100}[0,0]{$v_1^2$}[2,0]
	\mor{000}{010}[0,0]{$v_2^2$}[1,0]
	\mor{000}{001}[0,0]{$v_1^1$}[2,1]
	\mor{010}{110}[0,0]{$\kern-30pt v_1^2$}[2,0]
	\mor{100}{110}[0,0]{$v_2^2$}[1,0]
	\mor{100}{101}[0,0]{$v_1^1$}[2,1]
	\mor{010}{011}[0,0]{}[2,1]
	\mor{110}{111}[0,0]{$v_1^1$}[2,1]
	\mor{001}{011}[0,0]{$v_2^2$}[1,0]
	\mor{001}{101}[0,0]{$\kern30pt v_1^2$}[2,0]
	\mor{101}{111}[0,0]{$v_2^2$}[1,0]
	\mor{011}{111}[0,0]{$v_1^2$}[2,0]
	\enddc},
	\eqno(24)$$

\vskip15mm
\noindent  illustrates the composition
	$$\zd^3 f(q; v_1^1,v_1^2,v_2^2) = \zd^1\zd^2 f(q; v_1^1,v_1^2,v_2^2) = \zd^2 f(q + v_1^1;     v_1^2,v_2^2) - \zd^2 f(q;v_1^2,v_2^2). \eqno(25)$$

The relation
	$$\zd^n f(q;v_1, \dots,  v_n) = 
	\zd^{n-1} f(q + v_n;v_1, \dots,  v_{n-1}) - \zd^{n-1} f(q;v_1, \dots,  v_{n-1}) \eqno(26)$$
is a special case of Proposition 1.  The formula
		$$\zd^n f(q;v_1, \dots,  v_{n-1},0) = 0 \eqno(27)$$
is a consequence of this relation.

      \vskip 1\baselineskip 
	\proclaim Proposition 2.												
If $\zd^n f$ is multilinear at $q$, then $\zd^{n+1}f(q; \dots)=0$.

\noindent{\it Proof. }
	$$\eqalign{
	\zd^{n+1} f(q;w_1,w_2, v_1, \dots ,& v_{n-1}) \cr
	&= \zd^2\zd^{n-1} f(q;w_1,w_2, v_1, \dots , v_{n-1}) \cr
		&= \zd^{n-1} f(q+w_1+w_2; v_1, \dots , v_{n-1})  - \zd^{n-1} f(q+w_1; v_1, \dots , v_{n-1}) \cr
		&\phantom{=}  -\zd^{n-1} f(q+ w_2; v_1, \dots , v_{n-1}) + 
	\zd^{n-1} f(q; v_1, \dots , v_{n-1}) \cr  
	&= \zd^{n-1} f(q+w_1+w_2; v_1, \dots , v_{n-1})  - \zd^{n-1} f(q; v_1, \dots , v_{n-1}) \cr
	&\phantom{=} -\zd^{n-1} f(q+ w_1; v_1, \dots , v_{n-1}) + 
	 \zd^{n-1} f(q; v_1, \dots , v_{n-1}) \cr
	&\phantom{=} - \zd^{n-1} f(q+w_2;v_1, \dots , v_{n-1}) + 
	   \zd^{n-1} f(q; v_1, \dots , v_{n-1}) \cr
	&= \zd \zd^{n-1} f(q; w_1+w_2, v_1, \dots , v_{n-1}) - \zd\zd^{n-1} f(q; w_1, v_1, \dots , v_{n-1})\cr 
	&\phantom{=} - \zd\zd^{n-1} f(q; w_2, v_1, \dots , v_{n-1}) \cr
	&= \zd^n f (q; w_1+w_2,v_1, \dots , v_{n-1}) \cr
	&\phantom{=}  - \zd^n f(q; w_1, v_1, \dots , v_{n-1})  - \zd^n f(q; w_2, v_1, \dots , v_{n-1}) \cr
	&= 0 \hskip9.5cm\bull
        }\eqno(28)$$

   \vskip 1\baselineskip		
   \proclaim Proposition 3 (Leibniz rule).													
 The relation
       $$\zd^n(ff') (q; v_1, \dots, v_n) = \sum_{\scriptstyle I,I' \atop\scriptstyle I\cup I' = N}
             \zd^{|I|} f (q; {\b v}^{I})\zd^{|I'|} f' (q; {\b v}^{I'}) \eqno(29)$$
or
	$$\eqalign{
	\zd^n(ff') (q; v_1, \dots, v_n) 
	&= f(q)\zd^{n}f'(q;v_1,\dots, v_{n}) + \zd^{n}f(q;v_1,\dots, v_{n})f'(q) \cr
	&\hskip7mm + \sum_{\scriptstyle h,k =1,\dots, n
                     \atop\scriptstyle \{i_1,\dots, i_h\}\cup\{j_1,\dots, j_k\}=\{1,\dots, n\}}
           \zd^{h} f (q; v_{i_1}, \dots, v_{i_h})\zd^{k} f' (q; v_{j_1}, \dots, v_{j_k})
	}\eqno(30)$$
holds for any two functions $f \in A(Q)$ and $f' \in A(Q)$ and any $n \in \N$.

\noindent {\it Proof.} \enspace
	The formula is proved by induction.  
	In fact it is obviously valid for $n=0$.  For $n=1$ we have
	$$\eqalign{
	\zd(ff')(q;v) 
	&= ff'(q + v) - ff'(q) \cr
	&= f(q+v)f'(q+v) - f(q)f'(q) \cr
	&= \left(\zd f (q;v) + f(q) \right) ( \zd f'(q;v) + f'(q) ) - f(q)f'(q) \cr
	&= \zd f (q;v) \zd f'(q;v) +  \zd f(q;v) f'(q) + f(q) \zd f'(q;v)
	}\eqno(31)$$

	Assuming the formula valid for $n-1$, we have
		
	$$\eqalign{
	\zd^n(ff')
	&(q; v_1, \dots, v_n) = \zd\zd^{n-1}(ff') (q; v_1, \dots, v_n) \cr 
	&= \zd^{n-1}(ff') (q+v_n; v_1, \dots, v_{n-1}) - \zd^{n-1}(ff') (q; v_1, \dots, v_{n-1}) \cr
	&= f(q+v_n)\zd^{n-1}f'(q+v_n; v_1,\dots, v_{n-1})+ \zd^{n-1}f(q+v_n;v_1,\dots, v_{n  -1})f'(q+v_n)\cr   
	 &\hphantom{=} \cr
       &\qquad +\sum_{\scriptstyle h,k =1,\dots,n-1
         \atop\scriptstyle \{i_1,\dots, i_h\}\cup\{j_1,\dots, j_k\}=\{1,\dots, n-1\}}
         \zd^h f (q+v_n; v_{i_1}, \dots, v_{i_h}) \, \zd^k f '(q+v_n; v_{j_1}, \dots, v_{j_k}) \cr
	&\qquad -f(q)\zd^{n-1}f'(q;v_1,\dots, v_{n-1})-\zd^{n-1}f(q;v_1,\dots, v_{n-1})f'(q) \cr
	&\qquad  -\sum_{\scriptstyle h,k =1,\dots,n-1
         \atop\scriptstyle \{i_1,\dots, i_h\}\cup\{j_1,\dots, j_k\}=\{1,\dots, n-1\}} \zd^h f (q; v_{i_1}, \dots,
         v_{i_h}) \, \zd^k f '(q; v_{j_1}, \dots, v_{j_k}) \cr
         \cr
         &= (f(q)+\zd f(q;v_n))(\zd^{n}f'(q;v_1,\dots, v_{n})+\zd^{n-1}f'(q;v_1,\dots, v_{n-1})) \cr
       &\qquad +(\zd^{n}f(q;v_1,\dots, v_{n})+\zd^{n-1}f(q;v_1,\dots, v_{n-1}))(f'(q)+\zd f'(q;v_n)) \cr
       &\qquad +\sum_{\scriptstyle h,k =1,\dots,n-1
         \atop\scriptstyle \{i_1,\dots, i_h\}\cup\{j_1,\dots, j_k\}=\{1,\dots, n-1\}}
         \left[\zd^{h+1} f (q; v_{i_1}, \dots, v_{i_h}, v_n)+\zd^hf (q; v_{i_1}, \dots, v_{i_h})\right] \cr
       &\qquad\qquad\qquad\qquad\qquad\qquad
         \left[\zd^{k+1} f' (q; v_{j_1}, \dots, v_{j_k}, v_n)+\zd^{k}f '(q; v_{j_1}, \dots, v_{j_k})\right] \cr
       &\qquad -f(q)\zd^{n-1}f'(q,v_1,\dots, v_{n-1})-\zd^{n-1}f(q;v_1,\dots, v_{n-1})f'(q) \cr
       &\qquad -\sum_{\scriptstyle h,k =1,\dots,n-1 \atop\scriptstyle \{i_1,\dots, i_h\}
         \cup\{j_1,\dots, j_k\}=\{1,\dots, n-1\}} 
         \zd^h f (q; v_{i_1}, \dots, v_{i_h}) \, \zd^k f '(q; v_{j_1}, \dots, v_{j_k}) \cr
          }$$
	
	$$\eqalign{
       &=f(q)\zd^{n}f'(q;v_1,\dots, v_{n})+f(q)\zd^{n-1}f'(q;v_1,\dots, v_{n-1}) \cr
       &\qquad+\zd^{n}f(q;v_1,\dots, v_{n})f'(q)+\zd^{n-1}f(q;v_1,\dots, v_{n-1})f'(q) \cr
       &\qquad-f(q)\zd^{n-1}f'(q;v_1,\dots, v_{n-1})-\zd^{n-1}f(q;v_1,\dots, v_{n-1})f'(q) \cr
       &\qquad+\zd f(q;v_n)(\zd^{n}f'(q;v_1,\dots, v_{n})+\zd^{n-1}f'(q;v_1,\dots, v_{n-1})) \cr
	&\qquad+(\zd^{n}f(q;v_1,\dots, v_{n})+\zd^{n-1}f(q;v_1,\dots, v_{n-1}))\zd f'(q;v_n) \cr
       &\qquad +\sum_{\scriptstyle h,k =1,\dots,n-1
         \atop\scriptstyle \{i_1,\dots, i_h\}\cup\{j_1,\dots, j_k\}=\{1,\dots, n-1\}} \left\lbrack\zd^{h+1}
        f (q; v_{i_1}, \dots, v_{i_h}, v_n)\zd^{k+1} f' (q; v_{j_1}, \dots, v_{j_k}, v_n)\right. \cr
       &\qquad+\zd^{h} f (q; v_{i_1}, \dots, v_{i_h})\zd^{k+1} f' (q; v_{j_1}, \dots, v_{j_k}, v_n) \cr
       &\qquad+\zd^{h+1} f (q; v_{i_1}, \dots, v_{i_h}, v_n)\zd^{k} f' (q; v_{j_1}, \dots, v_{j_k}) \cr
       &\qquad\left.+\zd^{h} f (q; v_{i_1}, \dots, v_{i_h})\zd^{k} f' (q; v_{j_1}, \dots, v_{j_k})\right\rbrack \cr
       &\qquad-\sum_{\scriptstyle h,k =1,\dots,n-1 \atop\scriptstyle \{i_1,\dots, i_h\}
         \cup\{j_1,\dots, j_k\}=\{1,\dots, n-1\}} 
         \zd^h f (q; v_{i_1}, \dots, v_{i_h}) \, \zd^kf '(q; v_{j_1}, \dots, v_{j_k}) \cr
         \cr
       &=f(q)\zd^{n}f'(q;v_1,\dots, v_{n})+\zd^{n}f(q;v_1,\dots, v_{n})f'(q) \cr
       &\qquad+\zd f(q;v_n)\zd^{n}f'(q;v_1,\dots, v_{n})+\zd f(q;v_n)\zd^{n-1}f'(q;v_1,\dots, v_{n-1}) \cr
       &\qquad+\zd^{n}f(q;v_1,\dots, v_{n})\zd f'(q;v_n)+\zd^{n-1}f(q;v_1,\dots, v_{n-1})\zd f'(q;v_n) \cr
       &\qquad+\sum_{\scriptstyle h,k =1,\dots,n-1
         \atop\scriptstyle \{i_1,\dots, i_h\}\cup\{j_1,\dots, j_k\}=\{1,\dots, n-1\}}
       \left\lbrack\zd^{h+1} f (q; v_{i_1}, \dots, v_{i_h}, v_n)\zd^{k+1} f' (q; v_{j_1}, \dots, v_{j_k}, v_n)\right.\cr
        &\qquad+\zd^{h} f (q; v_{i_1}, \dots, v_{i_h})\zd^{k+1} f' (q; v_{j_1}, \dots, v_{j_k}, v_n)\cr
        &\qquad\left.+\zd^{h+1} f (q; v_{i_1}, \dots, v_{i_h}, v_n)
        \zd^{k} f' (q; v_{j_1}, \dots, v_{j_k})\right\rbrack \cr
        \cr
        &=f(q)\zd^{n}f'(q;v_1,\dots, v_{n})+\zd^{n}f(q;v_1,\dots, v_{n})f'(q)\cr
       &\qquad+\sum_{\scriptstyle h,k =1,\dots,n
         \atop\scriptstyle \{i_1,\dots, i_h\}\cup\{j_1,\dots, j_k\}=\{1,\dots, n\}}
         \zd^{h} f (q; v_{i_1}, \dots, v_{i_h})\zd^{k} f' (q; v_{j_1}, \dots, v_{j_k})
		}\eqno(32)$$
\hfill$\bull$		

\Sezione Polarization of mappings.

The polarization operators can be extended to mappings between affine spaces.  
The chain rule will be established for polarizations of compositions of mappings.

Let $g \colon P\to Q$ be a mapping from an affine space $(P,U)$ to an affine space $(Q,V)$.  
The $0$-{\it th polarization} of $g$ is the mapping itself.  
For $n>0$, the $n$-{\it th polarization} of $g$ is the mapping
	$$\zd^n g \colon P\times U^n \to V. \eqno(33)$$
defined by
	$$\eqalign{
	\zd^n g(p; u_1, \dots, u_n ) 
	&= (-1)^n\sum_{m=0}^n (-1)^m \sum_{\scriptstyle I \subset \{1,2,\ldots,n\}
	   \atop{\scriptstyle |I| = m}}
        g\left(p + {\s{\sum_{i \in I}}}u_i\right) \cr
	&= (-1)^n  \sum_{m=0}^n (-1)^m\; \sum_{1\le_1<\dots<i_m\le n} g(p + u_{i_1} +\dots + u_{i_m}).
	}\eqno(34)$$
	
The $n$-th unidirectional polarization
		$$\zD^n g \colon P \times U \to V \eqno(35)$$
of the mapping $g$ is defined by
		$$\zD^n g(p;u) = \zd^n g(p;u, \dots, u).\eqno(36)$$

\vskip1\baselineskip		
\proclaim Proposition 4 (The chain rule).													
Let $f$ be a function on an affine space $(Q,V)$ and let $g$ be a mapping from an affine space $(P,U)$ to $(Q,V)$.  For each $n \in \N$ we have the
	$$\zd^n(f \circ g)(p;u_1,\ldots,u_n) = \sum_{\scriptstyle  \{I^1,\ldots,I^k\}
										\atop{\scriptstyle {\bigcup}_{i=1}^k I^i\, = \, N
										\atop{\scriptstyle \,I^i \ne I^j {\rm if }\, i \ne j }}}
		 \zd^k f\left(g(p);\zd^{|I^1|}g(p; {\bu}^1),\ldots,\zd^{|I^k|}g(p; {\bu}^k)\right). \eqno(37)$$
The symbol ${\bu}^i$ stands for ${\bu}^{I^i}$.

\noindent {\it Proof.} \enspace
The formula
	$$\zd^n(f \circ g)(p;u_1,\ldots,u_n) = (-1)^n \sum_{I \subset N}(-1)^{|I|} f
	     \left(g\left(p +{\s{\sum_{i \in I}}}u_i\right)\right) \eqno(38)$$
follows from the definition of the $n$-th polarization.  
The expression
	$$g\left(p + {\s{\sum_{i \in I}}}u_i\right),  \eqno(39)$$
which is the argument of $f$ in the above formula, is transformed in
$${\s{\sum_{J \subset I}}}\zd^{|J|}g(p; {\bu}^J)  \eqno(40)$$
by applying formula (12).  The next step is the application of formula (12) to the expression
	$$f\left({\s{\sum_{J \subset I}}}\zd^{|J|}g(p; {\bu}^J)\right) \eqno(41)$$
with the result
	$$\sum_{\scriptstyle  \{I^1,\ldots,I^k\}
		\atop{\scriptstyle {\bigcup}_{i=1}^k I^i\, \subset \, I
		\atop{\scriptstyle \, I^i \ne I^j {\rm if }\, i \ne j }}}
		\zd^k f\left(g(p);\zd^{|I^1|}g(p; {\bu} ^1),\ldots,\zd^{|I^k|}g(p; {\bu}^k)\right)  \eqno(42)$$
The result of the transformations is the new version
	$$\zd^n(f \circ g)(p;u_1,\ldots,u_n) = (-1)^n \sum_{I \subset N}(-1)^{|I|} 
	     \sum_{\scriptstyle  \{I^1,\ldots,I^k\}
	    \atop{\scriptstyle {\bigcup}_{i=1}^k I^i\, \subset \, I
            \atop{\scriptstyle \, I^i \ne I^j {\rm if }\, i \ne j }}}
	     \zd^k f\left(g(p);\zd^{|I^1|}g(p; {\bu}^1),\ldots,\zd^{|I^k|}g(p; {\bu}^k)\right)  \eqno(43)$$
of formula (38).  The term
	$$\zd^k f\left(g(p);\zd^{|I^1|}g(p; {\it u}^1),\ldots,\zd^{|I^k|}g(p; {\bu}^k)\right)  \eqno(44)$$
in the above sum is repeated a number of times with different signs.  We are going to determine the coefficient of the combined term by simplifying the coefficient
	$$(-1)^n \sum_{\scriptstyle {\bigcup}_{i=1}^k I^i\, \subset \, I \subset N}(-1)^{|I|}.  \eqno(45)$$

Let $m$ denote the cardinality
	$$\left|{\bigcup}_{i=1}^k I^i\right|.  \eqno(46)$$

For each integer $h$ such that $m \le h \le n$, the number of sets $I$ of cardinality $|I| = h$ with
	$$\bigcup_{i=1}^k I^i\, \subset \, I \subset N  \eqno(47)$$
is equal to
	$${n-m \choose h-m}.  \eqno(48)$$
The equality
	$$\eqalign{
	(-1)^n \sum_{\scriptstyle {\bigcup}_{i=1}^k I^i\, \subset \, I \subset N}(-1)^{|I|}
	 &= (-1)^n\sum_{h=m}^n (-1)^h {n-m \choose h-m} \cr
	 &= \sum_{h=0}^{n-m} (-1)^{n-m-h} {n-m \choose h} \cr
	 \cr
	 &={ \cases {1 & {\rm if } $m=n$ \cr
			     0 & {\rm otherwise}}} 
        }\eqno(49)$$
implies that the value 1 for the coefficient (45) is obtained only when
	$$\bigcup_{i=1}^k I^i\, = \, I = N.  \eqno(50)$$
The resulting simplification of formula (43) proves the proposition.
\hfill$\bull$

\Sezione Polyhomogeneous functions.
		
\proclaim Definition 2.												
A function $f\in A(Q)$ is said to be {\it homogeneous of degree} $n$ at $q\in Q$ if
		$$f(q+\zl v) = \zl^n f(q+v)  \eqno(51)$$
for each $v \in V$ and $\zl \in {\R}$.

We denote by $Y^n(Q,q)$ the vector space of homogeneous functions of degree $n$ at $q$.  
If $f _1\in Y^{n_1}(Q,q)$ and $f _2\in Y^{n_2}(Q,q)$, then $f_1f_2 \in Y^{n_1+n_2}(Q,q)$. 

The formula
	$$\sum_{k=0}^n (-1)^{n-k} {n \choose k} (l + k)^m = 
	    \cases {\phantom{!} 0  & if $m < n$ \cr
						n!   & if $ m = n$\cr
	     }\eqno(52)$$
for integers $m$, $n$, and $l \ge 0$ was derived by Euler [Eu].  
The following relations are a consequence of (14) and formula (52) with $l = 0$. 
If $f\in Y^n(q,Q)$, then
	$$ \zD^j f(q;v)= 
		\cases { 0 &  if $ j > n $ \cr
                            j! f(q+v) &  if $ j = n $ \cr
              } \eqno(53)$$
and
       $$\zD^j f(q;\zl v) = \zl^n \zD^j f(q;v) \eqno(54)$$
for $j \le n$.  For $j = n$ the two relations imply
       $$\zD^n f(q;\zl v) = n! \zl^n f(q + v).  \eqno(55)$$

We will denote by $\zq^n$ the simple function
	$$\zq^n \colon {\R}\times{\R}\rightarrow{\R} \colon (q;\lambda) \mapsto (q+\lambda)^n. \eqno(56)$$
																							
\vskip1\baselineskip
\proclaim Proposition 5.												
If $f \in Y^n(Q,q)$, then
	$$\zD^k f(q+v;\zl v) = \zD^k\zq^n(1;\zl)f(q + v) \eqno(57)$$
for any integer $k$.
		
\noindent {\it Proof.}
	$$\eqalign{
		\zD^k f(q+v;\zl v) 
		&= (-1)^k \sum_{m=0}^k(-1)^m {k \choose m} f(q + v + m \zl v) \cr
		&= \left[ (-1)^k \sum_{m=0}^k(-1)^m {k \choose m}(1 + m \zl)^n \right] f(q + v) \cr
		&= \zD^k \zq^n(1; \zl) f(q + v).
	}\eqno(58)$$
The last equality follows from (14).
\hfill$\bull$

\vskip1\baselineskip
\proclaim Proposition 6.	
       $$\zD^k \zq^n(1;\zl) = \sum_{i=0}^n {n \choose i} \zl^i
	\left[(-1)^k \sum_{m=0}^k(-1)^m {k \choose m} m^i\right]. \eqno(59)$$
																								\noindent {\it Proof.}
	$$\eqalign{
	     \zD^k \zq^n(1;\zl) 
	     &= (-1)^k \sum_{m=0}^k(-1)^m {k \choose m}(1 + m \zl v)^n \cr
	     &= (-1)^k \sum_{m=0}^k(-1)^m {k \choose m}\sum_{i=0}^n {n \choose i} (m\zl)^i \cr
	     &= \sum_{i=0}^n {n \choose i} \zl^i\left[(-1)^k \sum_{m=0}^k(-1)^m {k \choose m} m^i\right].
	} \eqno(60)$$
	
\vskip1\baselineskip
\proclaim Corollary 1.																				
	$$\zD^k\zq^n(1;\zl) = \cases {0 &  if $ k > n $\cr
                                                     \zl^k k! &  if $ k = n $ 
                                                     } \eqno(61)$$

\noindent {\it Proof.} \enskip
If $n < k$, then the summation index $i$ in (59) is always less than $k$.  
It follows from the Euler formula (52) with $l = 0$ that all terms in the sum are zero.  
If $n = k$, then the only nonzero term in the sum is the one corresponding to $i = k$.  
This term is the expression $\zl^k k!$.
\hfill$\bull$

The relation
	$$\zD^k f(q+v;\zl v) = \sum_{i=0}^n {n \choose i}\zl^i
	\left[(-1)^k \sum_{m=0}^k(-1)^m {k \choose m} m^i\right] f(q + v). \eqno(62)$$
for $f \in Y^n(Q,q)$ is an immediate consequence of the last two propositions.  
The formula
	$$\zD^k f(q+v;\zl v) = \cases{ 0 &  if $ k > n $ \cr
                                                      \zl^k k! f(q + v) & if $ k = n $
                                                     }\eqno(63)$$
follows from Corollary 1.

\vskip1\baselineskip
\proclaim Definition 3. 											
{\it Polyhomogeneous functions} are sums of homogeneous functions.  
The {\it degree} of a polyhomogeneous function is the highest degree of its non zero homogeneous components.

The polarization operators are a tool for extracting homogeneous components of a polyhomogeneous function.  In fact, if
	$$f = \sum_{j=0}^n y^j \quad , \quad y^j\in Y^j(Q,q_0), \eqno(64)$$
then, owing to (53), for all $m\leq n$, 
	$$\eqalign{
		\zD^mf(q_0;q-q_0) 
		&= \sum_{j=0}^n \zD^m  y^j (q_0;q-q_0)\cr
		&= m! y^m \big(q_0 + (q-q_0)\big) + \zD^m \sum_{i=m+1}^n y^i (q_0;q-q_0).
		}\eqno(65)$$																				
Hence,
	$$y^m(q) = {1\over m!} \zD^{m}\left(f -  \sum_{j=m+1}^n y^j\right)(q_0;q-q_0). \eqno(66)$$
This formula extracts the $m$-th component after the components of higher degrees have been extracted.

\Sezione Polynomials.
		
\proclaim Definition 4.
A function $f\in A(Q)$ is said to be a {\it homogeneous polynomial of degree} $n$ at $q\in Q$ if there is a function $F\in B^n(Q)$, multilinear at $q$ such that, for each $v \in V$,
	$$f(q + v) =  {1\over n!} F(q; \underbrace{v,\dots, v}_{n\, {\rm times}}). \eqno(67)$$

We denote by $K^n(Q,q)$ the vector space of homogeneous polynomial functions of degree $n$ at $q$.  If $f _1\in K^{n_1}(Q,q)$ and $f _2\in K^{n_2}(Q,q)$, then $f_1f_2 \in K^{n_1+n_2}(Q,q)$. 

The function $F$ in the above definition is not unique.  A symmetric function can be chosen.  If $F$ satisfies relation (67), then its symmetric part satisfies the same relation.
  		
\vskip1\baselineskip		
\proclaim Proposition 7.
A function $f\in A(Q)$ is a homogeneous polynomial of degree $n$ at $q$ if and only if it is homogeneous of degree $n$ at $q$ and the polarization $\zd^n f$ is multilinear at $q$.

\noindent {\it Proof.} \enskip
(i) If $f$ is homogeneous and the polarization $\zd^n f$ is multilinear at $q$, then formula (53) implies
that $f$ is a homogeneous polynomial.

(ii) Conversely, if $f$ is a homogeneous polynomial of degree $n$ at $q$, then there is a function $F\in B^n(Q)$, multilinear at $q$ and symmetric, such that (67) holds.  It follows immediately that $f$ is homogeneous of degree $n$ at $q$, since
       $$f(q+\zl v) = {1\over n!}\, F(q; \underbrace{\zl v, \dots, \zl v}_{n \, \rm{ times}}) 
       = {1\over n!}\, \zl^n  F(q; \underbrace{v,\dots, v}_{n\,{\rm times}}) = \zl^n f(q+v). \eqno(68)$$

Let us consider the polarization
	$$\eqalign{
	\zd^n f (q; v_1, \dots, v_n) 
	&= \sum_{m=0}^n \; \sum_{1\le i_1<\dots<i_m\le n}(-1)^{n-m} f(q + v_{i_1} + \dots + v_{i_m})\cr
	&= {1\over n!}\,\sum_{m=0}^n \; \sum_{1\le i_1<\dots<i_m\le n}(-1)^{n-m} 
	     F(q ; v_{i_1} +.. + v_{i_m}, \dots, v_{i_1} +.. + v_{i_m}).
	}\eqno(69)$$																		
Owing to multilinearity and symmetry of $F$, we then have
	$$\eqalign{
	\zd^n f (&q; v_1, \dots, v_n) \cr
      &= \sum_{(i_1,\ldots,i_n) \in \{1,\ldots,m\}^n} F(q;v_{i_1},\ldots,v_{i_n}) \cr
	&= {1\over n!}\,\sum_{m=0}^n \; \sum_{1\le i_1<\dots<i_m\le n}(-1)^{n-m} 
	   \sum_{n_1+ \dots + n_m = n} {n! \over n_1! \dots n_m!} 
	  F(q ;\underbrace{v_{i_1} ,\dots, v_{i_1}}_{n_1 \,{\rm times}}, \dots, \underbrace{v_{i_m} , \dots,
         v_{i_m}}_{n_m \,{\rm times}})
	}\eqno(70)$$																				
The last passage consists in reordering the vectorial arguments of $F$ and counting the numbers of identical terms obtained.  The coefficient
	$${n! \over n_1!\dots n_m!} \eqno(71)$$
is the product
	$${n \choose n_1} {n - n_1 \choose n_2} {n - n_1 - n_2 \choose n_3}  \dots 
	     {n - n_1 - n_2 - \dots - n_{m-1} \choose n_m}. \eqno(72)$$

We can write
	$$\zd^n f (q; v_1, \dots, v_n) = \hat S_{v_1} + S_{v_1}, \eqno(73)$$	
where $\hat S_{v_1}$ is the sum of all the terms in (70) with $v_1$ absent and $S_{v_1}$ is the sum of all the terms in (70) containing $v_1$.  
If $v_1=0$, then $\zd^n f (q; v_1, \dots, v_n)=0$ by (27) and $S_{v_1}=0$. 
Hence $\hat S_{v_1}=0$ for all $(v_1, \dots, v_n)\in V^n$, and the equality
	$$\zd^n f (q; v_1, \dots, v_n) = S_{v_1} \eqno(74)$$	
is established.

Let $\hat S_{v_1, v_2}$ denote the part of $S_{v_1}$ with $v_2$ absent and $S_{v_1, v_2}$ be the sum of the remaining terms in $S_{v_1}$.  
Setting $v_2=0$, we have $\zd^n f (q; v_1, \dots, v_n)=0$ and $S_{v_1, v_2}=0$.  
Hence $\hat S_{v_1, v_2}=0$ and the new equality
	$$\zd^n f (q; v_1, \dots, v_n) = S_{v_1, v_2} \eqno(75)$$
is established.  Continuing the process we arrive at the equality
	$$\zd^n f (q; v_1, \dots, v_n) = S_{v_1, v_2, \dots, v_n}\;, \eqno(76)$$
where $S_{v_1, v_2, \dots, v_n}$ is the part of (73) with all the vectors $v_1, v_2, \dots, v_n$ present.

From
	$$S_{v_1, v_2, \dots, v_n} = F(q; v_1, \dots, v_n) \eqno(77)$$
it easily follows that the polarization
	$$\zd^n f (q; v_1, \dots, v_n) = F(q; v_1, \dots, v_n) \eqno(78)$$
is multilinear at $q$.
\hfill$\bull$

\vskip1\baselineskip
\proclaim Proposition 8.	
If $f \in K^n(Q,q)$, then
	$$f(q + v_1 + \dots + v_m) = 
	\sum_{\scriptstyle(n_1,\ldots,n_m) \in {\N}^m 
	\atop{\scriptstyle \! n_1 +\dots + n_m = n}}
	{1 \over n_1!\dots n_m!} \zd^n f(q;\underbrace{v_1, \dots, v_1}_{n_1\, {\rm times}} ,\ldots,     \underbrace{v_m,\dots, v_m}_{n_m \, {\rm times}}). \eqno(79)$$

\noindent {\it Proof.}\enskip
	$$\eqalign{
		n!f(q + v_1 + \dots + v_m) 
		&= \zD^n f(q;v_1 + \dots + v_m) \cr
		&= \zd^n f(q;v_1 + \dots + v_m, \ldots ,v_1 + \dots + v_m,) \cr
		&= \sum_{\scriptstyle(n_1,\ldots,n_m) \in {\N}^m 
		      \atop{\scriptstyle \! n_1 +\dots + n_m = n}}
                  {n! \over n_1!\dots n_m!} \zd^n 
 f(q;\underbrace{v_1, \dots, v_1}_{n_1\, {\rm times}},\ldots,\underbrace{v_m, \dots, v_m}_{n_m\, {\rm times}}).
		}\eqno(80)$$																				
	
We are using formula (53) and multilinearity and symmetry of the polarization $\zd^n f$.  
\hfill$\bull$

\vskip1\baselineskip
\proclaim Definition 5.
{\it Polynomials} are sums of homogeneous polynomials.  
The {\it degree} of a polynomial is the highest degree of its non zero homogeneous components.

The next result shows that polynomials do not need reference points for their definition.

\vskip1\baselineskip
\proclaim Proposition 9.
If $f\in A(Q)$ is a polynomial of degree $n$ at $q_0$, then it is a polynomial of degree $n$ at any other point $q$.

\noindent {\it Proof.}\enskip
Let $f$ be a polynomial of degree $n$ at $q_0$, i.e.,
	$$f=\sum_{m=0}^n f^m  \eqno(81)$$
	$$f^m(q_0+v) = {1\over m!}F^m(q_0;v,\dots,v)  \eqno(82)$$
with $F^m$ multilinear at $q_0$.

A direct computation shows that
	$$f=\sum_{m=0}^n g^m  \eqno(83)$$
	$$g^m(q+v) = {1\over m!}G^m(q;\underbrace{v,\dots,v}_{m\, {\rm times}}),  \eqno(84)$$
where
	$$G^m = \sum_{j=m}^n \, \sum_{1\le i_1<\dots<i_m\le j}\, F^j_{(i_1\dots i_m)}.  \eqno(85)$$
Here we consider the function
       $$ F^j_{(i_1\dots i_m)} : Q\times V^m \lra \R \eqno(86)$$
such that, for all $q\in Q$ and $v_1,\dots, v_m \in V$,
	$$F^j_{(i_1\dots i_m)}(q;v_1,\dots, v_m) = F^j(q_0; q-q_0,\dots,q-q_0, v_1, 
	q-q_0,\dots, q-q_0, v_m,q-q_0, \dots, q-q_0)  \eqno(87)$$
where, in the argument of the function on the right-hand side, the vectors
 $v_1,\dots, v_m $ occupy the positions with the indices $i_1,\dots, i_m$ respectively.
 \hfill$\bull$

We denote by $P^n(Q)$ the space of polynomials of degrees not  exceeding $n$.

The space 
	$$P(Q) = +_{n=0}^{\infty} P^n(Q) \eqno(88)$$
of all polynomials on $Q$ is a subalgebra of $A(Q)$.

Let $Q^n(Q)$ be the complement of $P^n(Q)$ in $P(Q)$; it is composed of polynomials whose homogeneous components are of degrees greater than $n$ and it is an ideal in $P(Q)$.  
The space $P^n(Q)$ is not a subalgebra, but it acquires the structure of an algebra through a natural identification with the quotient algebra $P(Q)/Q^n(Q)$. 

The space $P^n(Q)$ with this algebra structure is called the {\it algebra of truncated polynomials of degree} $n$.

\vskip 2\baselineskip
\goodbreak
\centerline {\bf{B.  Differential calculus in affine spaces.}}

Differential calculus in affine spaces will be presented as an application of the theory developed in part  A.

Starting from the definition of the $n$-th derivative as the infinitesimal limit of the $n$-th polarization, we will resume the preceding results and derive the classical properties of derivatives by taking the limits.

We will obtain formulae (such as Leibniz and chain rule) for multidirectional derivatives which seem to be more explicit and compact than those we know from literature (cf., e.g., [2], [6], [8]).

\Sezione Derivatives.

\proclaim Definition 6.
Let $f\in A(Q)$.  The limit
	$${\rd}^n f(q;v_1,v_2,\ldots,v_n) = \lim_{s\to 0}{1\over s^n}\zd^n f (q;sv_1,sv_2, \dots,sv_n),      
	 \eqno(89)$$
if it exists, is called the $n$-{\it th multidirectional derivative} of $f$ at the point $q\in Q$ in the
 {\it multidirection} $(v_1,v_2,\ldots,v_n) \in V^n$.

The $n$-{\it th directional derivative} at $q$ in the direction $v$ is the restriction
	$${\rD}^n f(q;v) = {\rd}^n f(q;v,v,\ldots,v)  \eqno(90)$$
of the multidirectional derivative to $Q$ times the diagonal of $V^n$.

Derivatives have the obvious homogeneity properties
	$${\rd}^n f(q;\zl v_1,\zl v_2,\ldots,\zl v_n) = \zl^n {\rd}^n f(q;v_1,v_2,\ldots,v_n)  \eqno(91)$$
and
	$${\rD}^n f(q;\zl v) = \zl^n {\rD}^n f(q;v).  \eqno(92)$$
	
The formula
	$${\rD}^n f(q;v) = \lim_{s\to 0} {1 \over s^n} \zD^n f(q;sv)
	                       = \lim_{s\to 0} {(-1)^n\over s^n}\, \sum_{m=0}^n (-1)^m {n\choose m} f(q+msv)
	 \eqno(93)$$
is a direct consequence of (14).  It will be used as an alternative version of the definition of the directional
derivative.

\vskip 1\baselineskip
\proclaim Theorem 1 (Leibniz rule).
If the derivatives ${\rd}^k f$ and ${\rd}^k f'$ of functions $f \in A(Q)$ and $f' \in A(Q)$ exist for $k \le n$, then
	$${\rd}^n(ff') (q; v_1, \dots, v_n) = \sum_{\scriptstyle I,I'
                     \atop{\scriptstyle I\cup I' = N \atop\scriptstyle \!\!I\cap I' = \emptyset}}
                 {{\rd}}^{|I|} f (q; {\rm{\bf v}}^{I}){\rd}^{|I'|} f' (q;{ \rm{\bf v}}^{I'})  \eqno(94)$$
or
	$$\eqalign{
		{\rd}^n(ff') (q; v_1, \dots, v_n) 
		&= f(q) {\rd}^{n}f'(q;v_1,\dots, v_{n}) + {\rd}^{n}f(q;v_1,\dots, v_{n})f'(q) \cr
		&\hskip7mm + \sum_{\scriptstyle h+k = n
                     \atop\scriptstyle \{i_1,\dots, i_h\}\cup\{j_1,\dots, j_k\}=\{1,\dots, n\}}
                {\rd}^{h} f (q; v_{i_1}, \dots, v_{i_h}) {\rd}^{k} f' (q; v_{j_1}, \dots, v_{j_k})
		}\eqno(95)$$																				
	
\noindent {\it Proof} \enskip
The proof is based on the polarization version of the Leibniz rule (Proposition 3).
	$$\eqalign{
	{\rd}^n&(ff') (q; v_1, \dots, v_n) = \lim_{s\to 0} {1\over s^n}\, \zd^n(ff') (q;sv_1, \dots, sv_n) \cr 
	&=\lim_{s\to 0} {1\over s^n}\, f(q)\zd^{n}f'(q;sv_1,\dots, sv_{n}) 
	+ \lim_{s\to 0} {1\over s^n}\, \zd^{n}f(q;sv_1,\dots, sv_{n})f'(q) \cr
	&\phantom{=} + \lim_{s\to 0} \!\!\!\!\sum_{\scriptstyle h,k =1,\dots,n
          \atop\scriptstyle \{i_1,\dots, i_h\}\cup\{j_1,\dots, j_k\}=\{1,\dots, n\}}
           s^{h+k-n}{1\over s^h} \zd^{h} f (q; sv_{i_1}, \dots, sv_{i_h}){1\over s^k} \zd^{k} f'(q; sv_{j_1},
\dots, sv_{j_k}) \cr
	&= f(q){\rd}^{n} f' (q; v_1, \dots, v_n) + {\rd}^{n} f (q; v_1, \dots, v_n)f'(q) \cr
	&\phantom{=} + \!\!\!\!\sum_{\scriptstyle h+k = n
         \atop\scriptstyle \{i_1,\dots, i_h\}\cup\{j_1,\dots, j_k\}=\{1,\dots, n\}}
         {{\rd}}^{h} f (q; v_{i_1}, \dots, v_{i_h}){\rd}^{k} f' (q; v_{j_1}, \dots, v_{j_k}).
	 }\eqno(96)$$	
We observe that in the first sum of the above formula the indices $h$ and $k$ satisfy the inequality $h + k \ge n$.
The terms with $h + k > n$ vanish since
	$$\lim_{s\to 0} s^{h+k-n} = 0. \eqno(97)$$	
If $h + k = n$, then $s^{h+k-n} = 1$.
\hfill$\bull$

\Sezione Derivatives of homogeneous functions.

\proclaim Proposition 10.
A function $f \in Y^n(Q,q)$ admits at $q$ directional derivatives of any order and in each direction and
	$${\rD}^m f(q;v) = m!\,\zd^m_n \,f(q + v),  \eqno(98)$$
where $\zd^k_n$ is the Kronecker symbol.

\noindent{\it Proof.}\enskip
If $m > n$, then
	$${\rD}^m f(q;v) = \lim_{s\to 0} {1 \over s^m} \zD^m f(q;sv) = 0  \eqno(99)$$
since,  by (53), $\zD^m f(q;sv) = 0$.  
If $m = n$, then
	$${\rD}^m f(q;v) = \lim_{s\to 0} {1 \over s^m} \zD^m f(q;sv)
	 = \lim_{s\to 0} {m! \over s^m} f(q + sv) = m!f(q + v)  \eqno(100)$$
since, by (53), $\zD^m f(q;sv) = m!f(q + sv) = m!s^mf(q + v)$.  
If $m < n$, then
	$${\rD}^m f(q;v) = \lim_{s\to 0} {1 \over s^m} \zD^m f(q;sv) 
	= \lim_{s\to 0} s^{n-m} \zD^m f(q;v) = 0  \eqno(101)$$
since $\zD^m f(q;sv) = s^n\zD^m f(q;v)$.
\hfill$\bull$

\vskip1\baselineskip
\proclaim Theorem 2 (Euler's formula).
If $f \in Y^n(Q,q)$, then
	$${\rD}^k f(q+v;v) = \cases{ 
	           0 & if $ k > n $\cr
	           \cr
                  \displaystyle{k!} f(q + v) &  if $k = n$ \cr
                  \cr
                  \displaystyle{ {n! \over (n-k)!}} f(q + v) &  if $k < n$ 
		}\eqno(102)$$
for each $v \in V$.
		
\noindent{\it Proof.}\enskip
The first two cases follow from (63).  
If $k \ge n$, then by (62) and (52),
	$$\eqalign{
	{\rD}^k f(q+v;v)& = \lim_{s\to 0} {1 \over s^k} \sum_{i=0}^n {n\choose i}s^i
	\left[(-1)^k\sum_{m=0}^k(-1)^m{k\choose m} m^i\right]f(q + v) \cr
	&= \lim_{s\to 0} \sum_{i=0}^{k-1}{n\choose i}{s^i\over s^k}
	\left[(-1)^k \sum_{m=0}^k(-1)^m{k\choose m} m^i\right]f(q + v) \cr
	&\hskip10mm + \lim_{s\to 0} {n\choose k} 
	\left[(-1)^k \sum_{m=0}^k(-1)^m{k\choose m} m^k\right]f(q + v) \cr
	&\hskip10mm + \lim_{s\to 0} \sum_{i=k+1}^{n}{n\choose i}{s^i\over s^k}
	\left[(-1)^k \sum_{m=0}^k(-1)^m {k\choose m} m^i\right]f(q + v) \cr
	&= {n\choose k} k!f(q + v) \cr
	&= {n! \over (n - k)!}f(q + v) .
	}\eqno(103)$$
\hfill$\bull$

\vskip1\baselineskip
\proclaim Corollary 2.
If $f \in Y^n(Q,q)$, then
	$${\rD}^k f(q+v;v) = {(n+m-k)! \over (n-k)!} {\rD}^{k-m}f(q + v;v) \eqno(104)$$
for $m \le k \le n$ and each $v \in V$.
\hfill$\bull$

Results derived for homogeneous functions are valid for homogeneous polynomials.  
For a homogeneous polynomial $f \in K^n(Q,q)$, in view of Proposition 7, formula (89) implies the formulae
	$${\rd}^n f(q;v_1,v_2,\ldots,v_n) = \zd^n f (q;v_1,v_2, \dots,v_n) \eqno(105)$$
and by (53)
	$${\rD}^n f(q;v) = \zD^n f(q;v) = n!f(q + v). \eqno(106)$$

\Sezione Iterated derivatives.
	
The formula
	$${\rd}^n F(q; v_1, \dots, v_n , w_1, \dots, w_m) 
	= \lim_{s\rightarrow  0} {1\over s^n} \zd^n F(q; sv_1, \dots, sv_n ,w_1, \dots, w_m) \eqno(107)$$
extends the definition (89) to a function $F \in B^m(Q)$ with the polarization defined in (17).  
This extension permits the iteration of the construction of derivatives.

Let us explicitly remark that  the existence of the $n$-th multidirectional derivative defined by (89) does not imply the existence of the derivatives of lower orders and, even in the case where these exist, they need not give the $n$-th derivative by iteration. The required additional condition is the
continuity of all the derivatives involved. Under this assumption we have the following results.

\vskip1\baselineskip
\proclaim Proposition 11.
The relation
	$${\rd}^{n_1}\,{\rd}^{n_2} = {\rd}^{n_1+n_2} \eqno(108)$$
holds for all integers $n_1$ and $n_2$ such that ${\rd}^{n_1}\,{\rd}^{n_2}f$ exists.

\noindent{\it Proof.}\enskip
       $$\eqalign{
	{\rd}^{n_1}{\rd}^{n_2} 
	&f(q; v^1_1, \dots , v^1_{n_1}, v^2_1, \dots , v^2_{n_2}) 
	= \lim_{s \rightarrow 0}{1\over s^{n_1}}\zd^{n_1}{\rd}^{n_2}f(q; sv^1_1, \dots , sv^1_{n_1}, v^2_1, \dots , v^2_{n_2})\cr
	&= \lim_{s_1 \rightarrow 0} \lim_{s_2 \rightarrow 0}{1\over s^{n_1}}{1\over s^{n_2}} \zd^{n_1}\zd^{n_2}
	f(q; s_1v^1_1, \dots ,s_1v^1_{n_1},s_2v^2_1, \dots , s_2v^2_{n_2}) \cr
	&= \lim_{s \rightarrow 0}{1\over s^{n_1+n_2}} \zd^{n_1+n_2}f(q; sv^1_1, \dots , sv^1_{n_1}, sv^2_1, \dots , sv^2_{n_2}) \cr
	&= {\rd}^{n_1+n_2} f(q; v^1_1, \dots , v^1_{n_1}, v^2_1, \dots , v^2_{n_2})
	}\eqno(109)$$
\hfill$\bull$

The following proposition establishes a useful relation between derivatives of functions on an affine space and derivatives of functions of real variables.
		
\vskip1\baselineskip
\proclaim Proposition 12.
If a function $f \in A(Q)$ admits the $n$-th multidirectional derivative in the multidirection $(v_1,v_2,\ldots,v_n) \in V^n$ at a point $q \in Q$, then
	$${\rd}^n f(q;v_1,v_2,\ldots,v_n) 
	= \left.{\partial^n\over\partial s_1\partial s_2\dots\partial s_n} 
	f(q + s_1 v_1 + s_2 v_2 + \dots + s_n v_n)\right|_{s_1=0,s_2=0,\ldots,s_n=0}. \eqno(110)$$																								
\noindent{\it Proof.}
	$$\eqalign{
	{\partial^n\over\partial s_1\partial s_2\dots\partial s_n} 
	&\left.\vphantom{{\partial\over\partial}} f(q + s_1 v_1 + s_2 v_2 + \dots + s_n v_n)
        \right|_{s_1=0,s_2=0,\ldots,s_n=0} \cr
	&= {\partial^{n-1}\over\partial s_2\dots\partial s_n}
	 \left.\lim_{s_1\rightarrow 0}{1\over s_1} \zd^1 f(q + s_2 v_2 + \dots + s_n v_n;s_1 v_1)
	 \right|_{s_2=0,\ldots,s_n=0} \cr
	&\hskip20mm ................................................. \cr
	&= \lim_{s_1\rightarrow 0}\lim_{s_2\rightarrow 0}\dots \lim_{s_n\rightarrow 0}
	  {1\over s_1 s_2 \dots s_n}\zd^1\zd^1\dots\zd^1 f(q;s_1 v_1,s_2 v_2,\dots,s_n v_n) \cr
	&= \lim_{s_1\rightarrow 0}\lim_{s_2\rightarrow 0}\dots \lim_{s_n\rightarrow 0}
	 {1\over s_1 s_2 \dots s_n} \zd^n f(q;s_1 v_1,s_2 v_2,\dots,s_n v_n) \cr
	&= \lim_{s_1\rightarrow 0}\lim_{s_2\rightarrow 0}\dots\lim_{s_n\rightarrow 0}
	\left[ {1\over s_1 s_2 \dots s_n} \zd^n f(q;s_1 v_1,s_2 v_2,\dots,s_n v_n)\right] \cr
	&= \lim_{s\rightarrow 0}\left[{1\over s^n} \zd^n f(q;sv_1,sv_2,\dots,sv_n)\right] \cr
	&= \lim_{s\rightarrow 0}{1\over s^n} \zd^n f(q;sv_1,sv_2,\dots,sv_n) \cr
	&= {\rd}^n f(q;v_1,v_2,\ldots,v_n).
		}\eqno(111)$$

We have replaced the expression
	$${1\over s_1 s_2 \dots s_n} \zd^n f(q;s_1 v_1,s_2 v_2,\dots,s_n v_n) \eqno(112)$$
by its continuous extension
	$$\eqalign{
	  &\left[ {1\over s_1 s_2 \dots s_n} \zd^n f(q;s_1 v_1,s_2 v_2,\dots,s_n v_n)\right] \cr
	  &\hskip20mm = \lim_{(s'_1,s'_2,\ldots,s'_n) \rightarrow (s_1,s_2,\ldots,s_n)}\quad
	  {1\over s'_1 s'_2 \dots s'_n}\zd^n f(q;s'_1 v_1,s'_2 v_2,\dots,s'_n v_n).
		}\eqno(113)$$																			
Continuity of this extension as a function of $(s_1,s_2,\ldots,s_n)$ is used.
\hfill$\bull$

\Sezione The chain rule.

\proclaim Theorem 3 (The chain rule).
Let $(P,U), (Q,V)$ be affine spaces and let  $g\colon P \rightarrow  Q$ and $f \in A(Q)$ be continuous mappings with continuous multidirectional derivatives of any order $k \le n$.  Then, with the notation of (36), for each $p \in P$ and $(u_1, \ldots , u_n) \in U^n$ we have
	$${\rd}^n(f\circ g)(p; u_1, \dots , u_n) 
	      = \sum_{\scriptstyle  \{I^1,\ldots,I^k\}
		\atop{\scriptstyle {\bigcup}_{i=1}^k I^i\, = \, N
		\atop{\scriptstyle  {\sum}_{i=1}^k | I^i |\, = \, n}}}
		{\rd}^k f\left(g(p); {\rd}^{|I^1|}g(p;{\bu}^1),\ldots, {\rd}^{|I^k|}g(p;{\bu}^k)\right). \eqno(114)$$

\noindent{\it Proof.}\enskip  
The equality
	$$\eqalign{
	{\rd}^n(f \circ g)(p; u_1, \dots , u_n) 
	&= \lim_{s\to 0}{1\over s^n} \sum_{\scriptstyle \{I^1,\ldots,I^k\}
			\atop{\scriptstyle {\bigcup}_{i=1}^k I^i\, = \, N
			\atop{\scriptstyle I^i \ne I^j {\rm if }\ i \ne j }}}
	       \zd ^k f\left(g(p); \zd^{|I^1|}g(p; s\bu^1),\ldots, \zd^{|I^k|}g(p; s\bu^k)\right) \cr
	&= \sum_{\scriptstyle \{I^1,\ldots,I^k\}
	  \atop{\scriptstyle {\bigcup}_{i=1}^k I^i\, = \, N
	  \atop{\scriptstyle I^i \ne I^j {\rm if }\ i \ne j }}}
	  \lim_{s\to 0}{1\over s^n}\zd ^k f\left(g(p); \zd^{|I^1|}f(q; s\bu^1),\ldots, \zd^{|I^k|}g(p; s\bu^k)\right)
		}\eqno(115)$$																			
follows from Proposition 4.  
Each individual term in the sum is computed using the continuous extensions introduced in the proof of Proposition 12.  
The result
	$$\eqalign{
	\lim_{s\to 0}{1\over s^n}\zd ^k
	& f\left(g(p); \zd^{m_1}f(q; s\bu^1),\ldots, \zd^{m_k}g(p; s\bu^k)\right) \cr
	&= \lim_{s\to 0}{1\over s^n}\zd ^k f\left(g(p); s^{m_1}
	  \left[{1\over s^{m_1}}\zd^{m_1}g(p; s\bu^1)\right],\ldots, s^{m_k}
         \left[{1\over s^{m_k}}\zd^{m_k}g(p; s\bu^k)\right]\right) \cr
	&= \lim_{s\to 0}\left[{1\over s^n}\zd ^k f\left(g(p); s^{m_1}
	  \left[{1\over s^{m_1}}\zd^{m_1}g(p; s\bu^1)\right],\ldots, s^{m_k}
         \left[{1\over s^{m_k}}\zd^{m_k}g(p; s\bu^k)\right]\right)\right] \cr
	&= \lim_{s\to 0}\lim_{t\to 0}\left[{1\over s^n}\zd ^k f\left(g(p);
          s^{m_1}\left[{1\over t^{m_1}}\zd^{m_1}g(p; t\bu^1)\right],\ldots, s^{m_k}
          \left[{1\over t^{m_k}}\zd^{m_k}g(p; t\bu^k)\right]\right)\right] \cr
	&= \lim_{s\to 0}\left[{1\over s^n}\zd ^k f
	  \left(g(p); s^{m_1}{\rd}^{m_1}g(p; \bu^1),\ldots, s^{m_k}{\rd}^{m_k}g(p; \bu^k)\right)\right] \cr
	&= \lim_{s\to 0}\left[{s^{m-n}\over s^{m_1}\dots s^{m_k}}\zd ^k f
	  \left(g(p); s^{m_1}{\rd}^{m_1}g(p; \bu^1),\ldots, s^{m_k}{\rd}^{m_k}g(p; \bu^k)\right)\right] \cr
	&= \lim_{\scriptstyle s \to 0
		\atop{\scriptstyle ( t_1 , \dots, t_k) \to (0, \dots, 0)}}\left[{s^{m-n}\over t_1\dots t_k}\zd ^k f
          \left(g(p); t_1{\rd}^{m_1}g(p; \bu^1),\ldots, t_k{\rd}^{m_k}g(p; \bu^k)\right)\right] \cr
	&= \lim_{\scriptstyle s \to 0
		\atop{\scriptstyle t \to 0}}\left[{s^{m-n}\over t^k}\zd ^k f
          \left(g(p); t{\rd}^{m_1}g(p; \bu^1),\ldots, t{\rd}^{m_k}g(p;\bu^k)\right)\right] \cr
	&= \lim_{s \to 0}s^{m-n}{\rd} ^k f
	\left(g(p);{\rd}^{m_1}g(p;\bu^1),\ldots, {\rd}^{m_k}g(p; \bu^k)\right) \cr
	&= \cases{ {\rd} ^k f\left(g(p);{\rd}^{m_1}g(p;\bu^1),\ldots, {\rd}^{m_k}g(p; \bu^k)\right) 
	                  &  if $m = n$\cr
			 0 &  if $m > n$ }
		}\eqno(116)$$																			with $m_i = |I^i|$ for $i = 1,\ldots,k$ and $m = m_1 + \dots + m_k$ inserted in (115) produces the proof.
\hfill$\bull$

\Sezione Differentiability.

The differentiability of a function, at any order, is related to the possibility of expressing it as a sum of a polynomial and a remainder.

We have already introduced the notion of polynomial function on an affine space. 
Now we will extend the notion of remainder. 
In the following definition, we will refer to one of the equivalent norms present in a vector space of finite dimension.

\vskip1\baselineskip
\proclaim Definition 7.
A function $r\in A(Q)$ is said to be a {\it remainder of order} $n$ at $q$ if $r(q)=0$ and
$$\lim_{v\to 0} {r(q+v)\over \| v \|^n}=0. \eqno(117)$$

We observe that a remainder of order $n$ at $q$ is also a remainder of any order $m<n$. This is a consequence of
$$\lim_{v\to 0} {r(q+v)\over \| v \|^m} = 
\lim_{v\to 0} \| v \|^{n-m} {r(q+v)\over \| v \|^n} = 0. \eqno(118)$$

\vskip1\baselineskip
\proclaim Proposition 13.
For any $k\le n$, we have
$${\rD}^kr(q;v) = 0 \eqno(119)$$

{\it Proof.} According to (93)
$${\rD}^kr(q;v) = \lim_{s \to 0} {\zD^kr(q;sv)\over s^k} =
 (-1)^k\sum_{l=0}^k(-1)^l {k\choose l}\lim_{s \to 0} {r(q+lsv)\over s^k}\,. \eqno(120)$$
But (117) implies
$$\lim_{s \to 0} {r(q+lsv)\over s^k} = 
 \| lv \|^n \lim_{s \to 0} s^{n-k}{r(q+lsv)\over \| lsv \|^n} = 0. \eqno(121)$$
 \hfill$\bull$

\vskip1\baselineskip
\proclaim Definition 8.
A function $f\in A(Q)$ is said to be  {\it differentiable} at $q\in Q$ if there is a polynomial $p\in P^1(Q)$ such that
$$r = f - p  \eqno(122)$$
is a remainder of order $1$ at $q$.

Via polarization, it follows from (122) that the homogeneous components of $p$ are
$$k^0 = f(q)  \eqno(123)$$
and
$$k^1(q+v) = {\rd}^1f(q;v)  \eqno(124)$$
The concept of differentiability extends easily to mappings between affine spaces.
For elements of the space $B^m(Q)$ we have the following definition.

\vskip1\baselineskip
\proclaim Definition 9.
A function $F \in B^m(Q)$ is said to be  {\it differentiable} at $q\in Q$ if there is a function $p\in B^m(Q)$  such that, for each $(w_1,\dots,w_m)\in V^m$,
$p(\cdot,w_1,\dots,w_m)\in P^1(Q)$ and
$$r(\cdot,w_1,\dots,w_m) = F(\cdot,w_1,\dots,w_m) -p(\cdot,w_1,\dots,w_m) \eqno(125)$$
is a remainder of order $1$ at $q$.

If a function $f$ is differentiable in $Q$, or at least in a neighbourhood  of a point $q$, the differentiability of the derivative ${\rd}^1f \in B^1(Q)$ can be examined.
A function $f\in A(Q)$ will be said to be twice differentiable at $q\in Q$ if it is differentiable in a neighbourhood of $q$ and the derivative ${\rd}^1f$ is differentiable at $q$. In such a case 
${\rd^1 \rd^1}f$ exists at $q$ and, owing to Proposition 11, it equals the second derivative ${\rd}^2f$ at the same point $q$.

Higher order differentiability will then be defined by induction.

\vskip1\baselineskip
\proclaim Definition 10.
A function $f \in A(Q)$ is said to be $n$  {\it times differentiable} at $q\in Q$ if it is $n-1$ times differentiable in a neighbourhood of $q$ and the derivative ${\rd}^{n-1}f$ is differentiable at $q$.

Differentiability of order $n$ assures the existence of all iterated derivatives up to order $n$.

\vskip1\baselineskip
The above definitions allow us to extend Taylor's theorems from  Banach spaces to affine spaces. For the sake of completeness we recall these results, whose claims we adapt to our context,  and whose proves can be found in standard textbooks of analysis (cf., e.g., [1],[2],[3],[4],[7],[9]).

\vskip1\baselineskip
\proclaim Theorem 4 (Taylor's formula).
If $f\in A(Q)$ is $n$ times differentiable at $q$, then there is a polynomial $p\in P^n(Q)$ such that
$$r = f - p  \eqno(126)$$
is a remainder of order $n$ at $q$.
\hfill$\bull$

The Taylor formula is the unique decomposition of a differentiable function as a sum of a polynomial and a remainder. In fact

\vskip1\baselineskip
\proclaim Proposition 14.
Let
$$f = \sum_{m=0}^n k^m + r \,,\eqno(127)$$
where $k^m\in K^m(Q,q)$ and $r $ is a remainder of order $n$ at $q$.
Then, for $0\le j\le n$,
$$k^j(q+v) = {1\over j!}{\rD}^j f(q ;v) \eqno(128)$$

{\it Proof.} Owing to (98) and (119), we have
$${\rD}^j f(q;v) = \sum_{m=0}^n {\rD}^j k^m(q;v) + {\rD}^j r(q;v)=j! k^j(q+v)\,.\eqno(129)$$
\hfill$\bull$

Then Taylor's formula assumes the form
$$f(q+v) = \sum_{m=0}^n{1\over m!}  {\rD}^m  f(q ;v) + R^n f(q;v) \eqno(130)$$
with
$$R^n f(q;v) = r(q+v) \eqno(131)$$

 Theorem 4 does not have a converse: the existence of such a decomposition does not imply differentiability. A converse theorem, however,  does  exist for a modified Taylor's theorem which concerns functions of class $C^n$.

\vskip1\baselineskip
\proclaim Definition 11.
A function $f \in A(Q)$ is said to be  {\it of class}  $C^n$  if it is $n$ times differentiable and its $n$-th derivative ${\rD}^n f : Q \times V \to \R$ is continuous.

\vskip1\baselineskip
\proclaim Theorem 5 (modified Taylor's theorem). 
If $f\in A(Q)$ is of class $C^n$, then for each $(q,q')\in Q\times Q$
$$f(q') = p(q,q') + r(q,q') ,\eqno(132)$$
where $p$ is the function
$$p : Q\times Q \lra {\R} ; (q,q')\mapsto \sum_{k=0}^n{1\over k!}{\rD}^k f(q;q'-q)\eqno(133)$$
and $r$ is a function on $Q\times Q$ such that
$$r(q,q)= 0 \quad {\sl and} \quad 
 \lim_{(q,q')\to(q_0,q_0)} {r(q,q') \over \|q'-q\|^n} = 0\eqno(134)$$
 for each $q_0$.
 \hfill$\bull$
 
 \vskip1\baselineskip
\proclaim Theorem 6 (converse of the modified Taylor's theorem). 
Let $p$ be a continuous function on $Q\times Q$ such that  for each $q$ the function $p(q,\cdot)$ is a polyomial of degree $n$ and let $r$ be function on $Q\times Q$ such that 
$$r(q,q)= 0 \quad {\sl and} \quad 
 \lim_{(q,q')\to(q_0,q_0)} {r(q,q') \over \|q'-q\|^n} = 0\eqno(135)$$
 for each $q_0$.
The function $f\in A(Q)$ defined by
$$f(q') = p(q,q') + r(q,q') \eqno(136)$$
is of class $C^n$ and
$$p : Q\times Q \lra {\R} ; (q,q')\mapsto p(q,q') =  \sum_{k=0}^n{1\over k!}D^k f(q;q'-q).\eqno(137)$$

\References

\noindent
\vskip1\baselineskip
$[1]$ Abraham, R., Marsden, J.E. and Ratiu, T., 1988, {\it Manifolds, tensor analysis and applications. Second edition.} (New York: Springer Verlag).
\vskip1\baselineskip
$[2]$ Abraham, R. and Robbin, J., 1967, {\it Transversal mappings and flows}  (New York--Amsterdam: Benjamin).
\vskip1\baselineskip
$[3]$ Albrecht, F. and Diamond, H.G., 1971, A converse of Taylor's theorem. {\it Indiana University Mathematics Journal}, {\bf 21}, 347--350.
\vskip1\baselineskip
$[4]$ Dieudonn\'e, J., 1960, {\it Foundations of modern analysis}  (New York: Academic Press).
\vskip1\baselineskip
$[5]$ Euler, L., 1754/55, Demonstratio theorematis Fermatiani omnem numerum primum formae 4n+1 esse summam duorum quadratorum. {\it Novi commentarii academiae scientiarum Petropolitanae}, {\bf 5}, 3--13. 
\vskip1\baselineskip
$[6]$ Fa\`a di Bruno, F., 1869, {\it Traite elementaire de calcul} (Paris: Gauthier--Villars).
\vskip1\baselineskip
$[7]$ Glaeser, G., 1958, Etude de quelques algebres tayloriennes. {\it Journal d'Analyse Math\'ematique}, {\bf 6}, 1--124.
\vskip1\baselineskip
$[8]$ Gradshteyn, I.S. and Ryzhik, I.M., 1965, {\it Tables of integrals, series and products}  (New York: Academic Press).
\vskip1\baselineskip
$[9]$ Nelson, E., 1969, {\it Topics in dynamics--I: Flows}  (Princeton: Princeton University Press).

\bye